\documentclass[12pt]{amsart}
\usepackage{amscd,amsmath,amssymb,amsfonts}

\usepackage{epsf}
\def\zigzag{\raisebox{-2cm}{\epsfysize=4cm\epsfbox{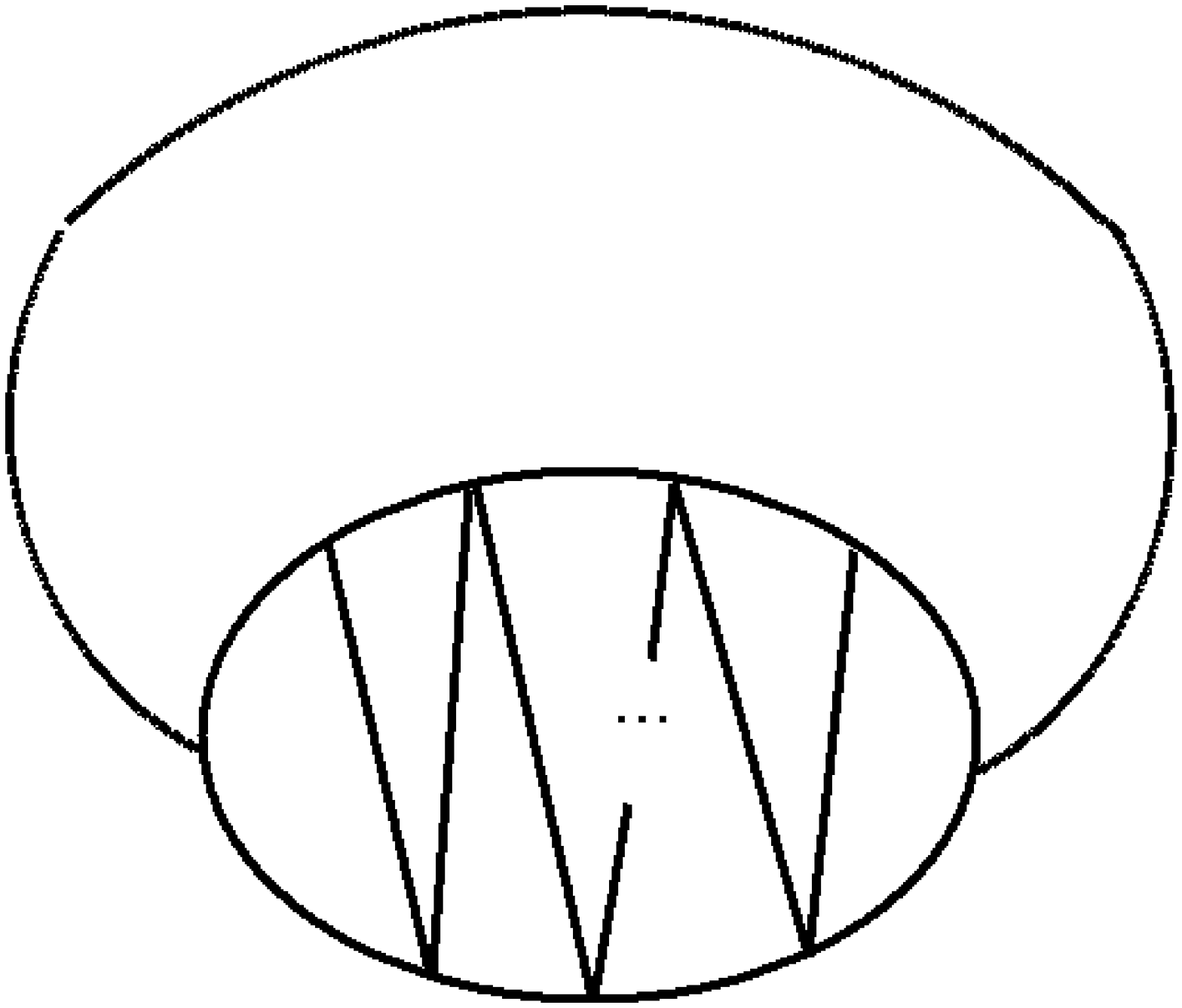}}\;}
\def\wheels{\raisebox{-2cm}{\epsfysize=4cm\epsfbox{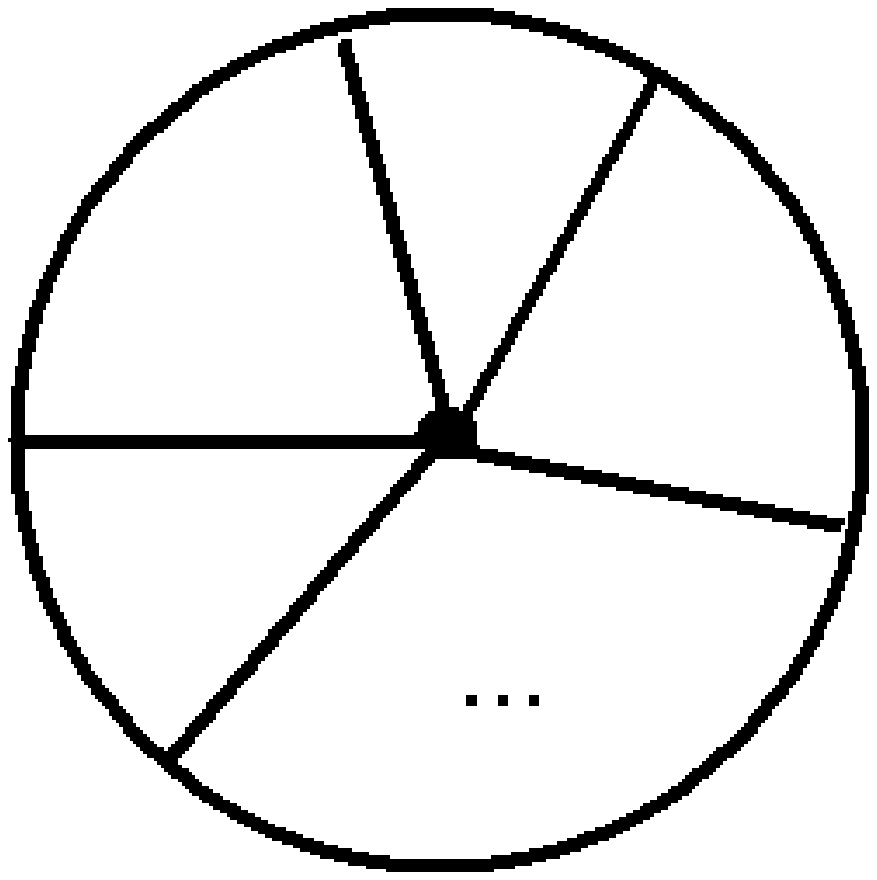}}\;}
\def\wthree{\raisebox{-2cm}{\epsfysize=3cm\epsfbox{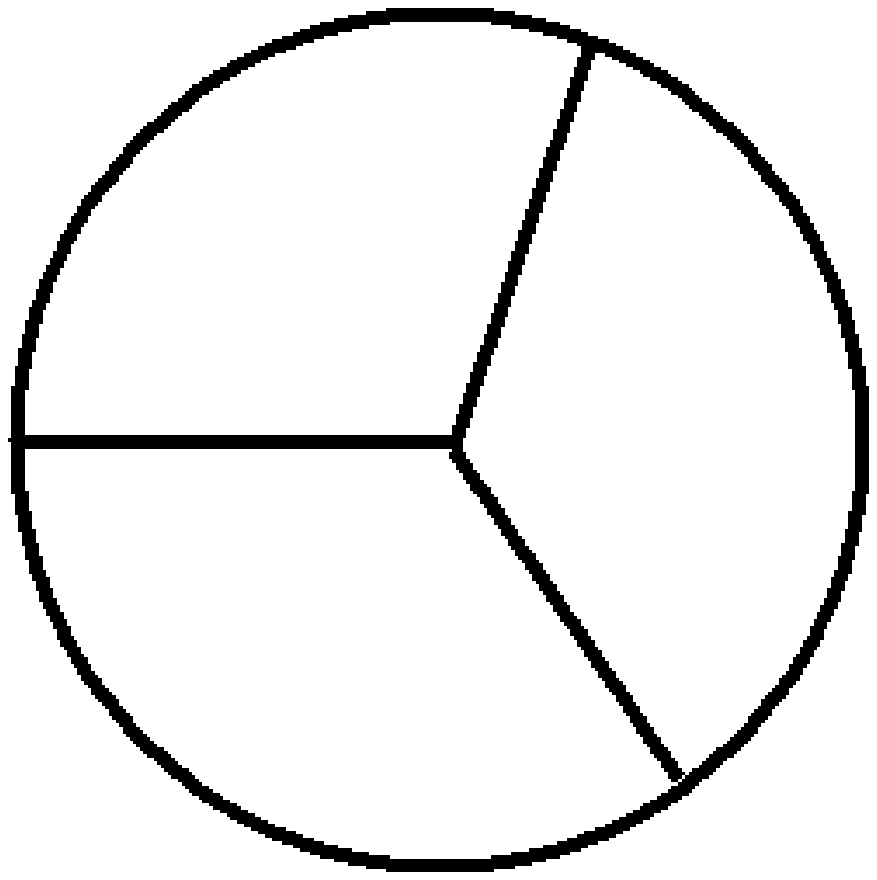}}\;}
\def\wfour{\raisebox{-2cm}{\epsfysize=3cm\epsfbox{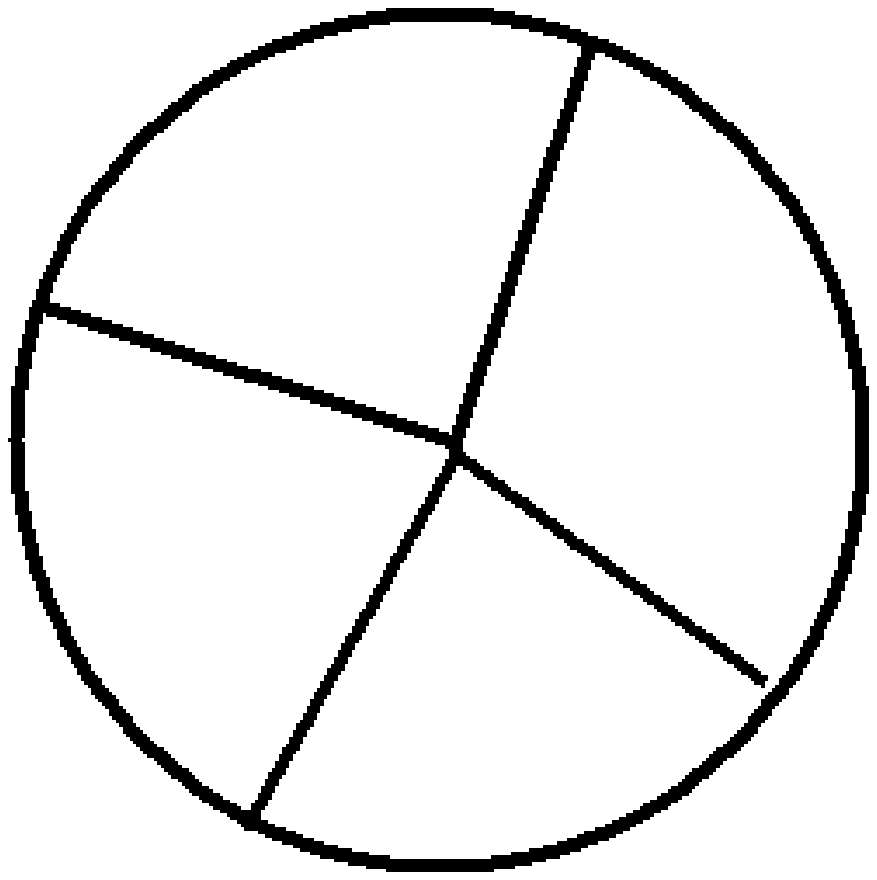}}\;}
\def\Mgraph{\raisebox{-2cm}{\epsfysize=4cm\epsfbox{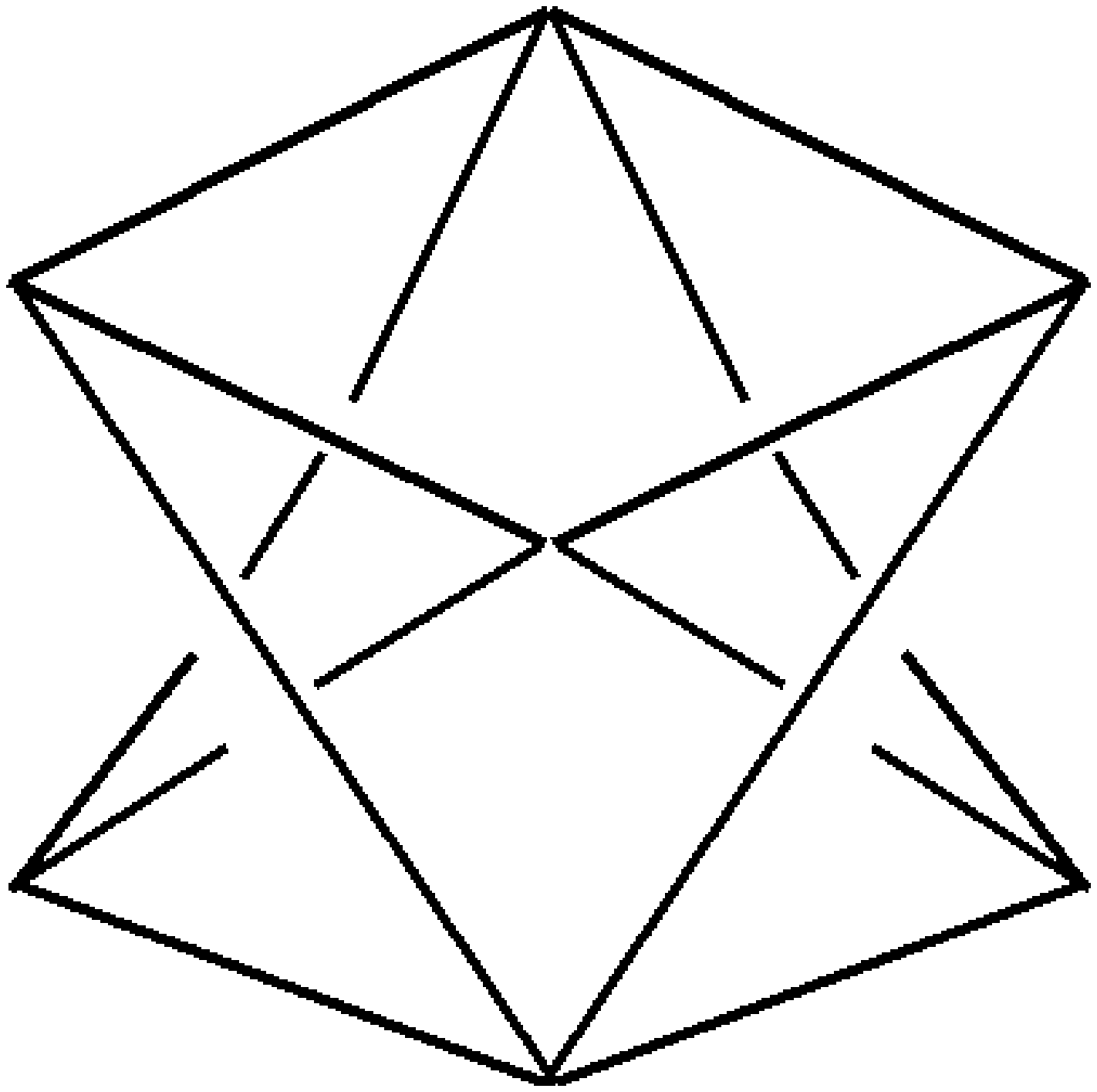}}\;}
\def\Migraphs{\raisebox{-2cm}{\epsfysize=4cm\epsfbox{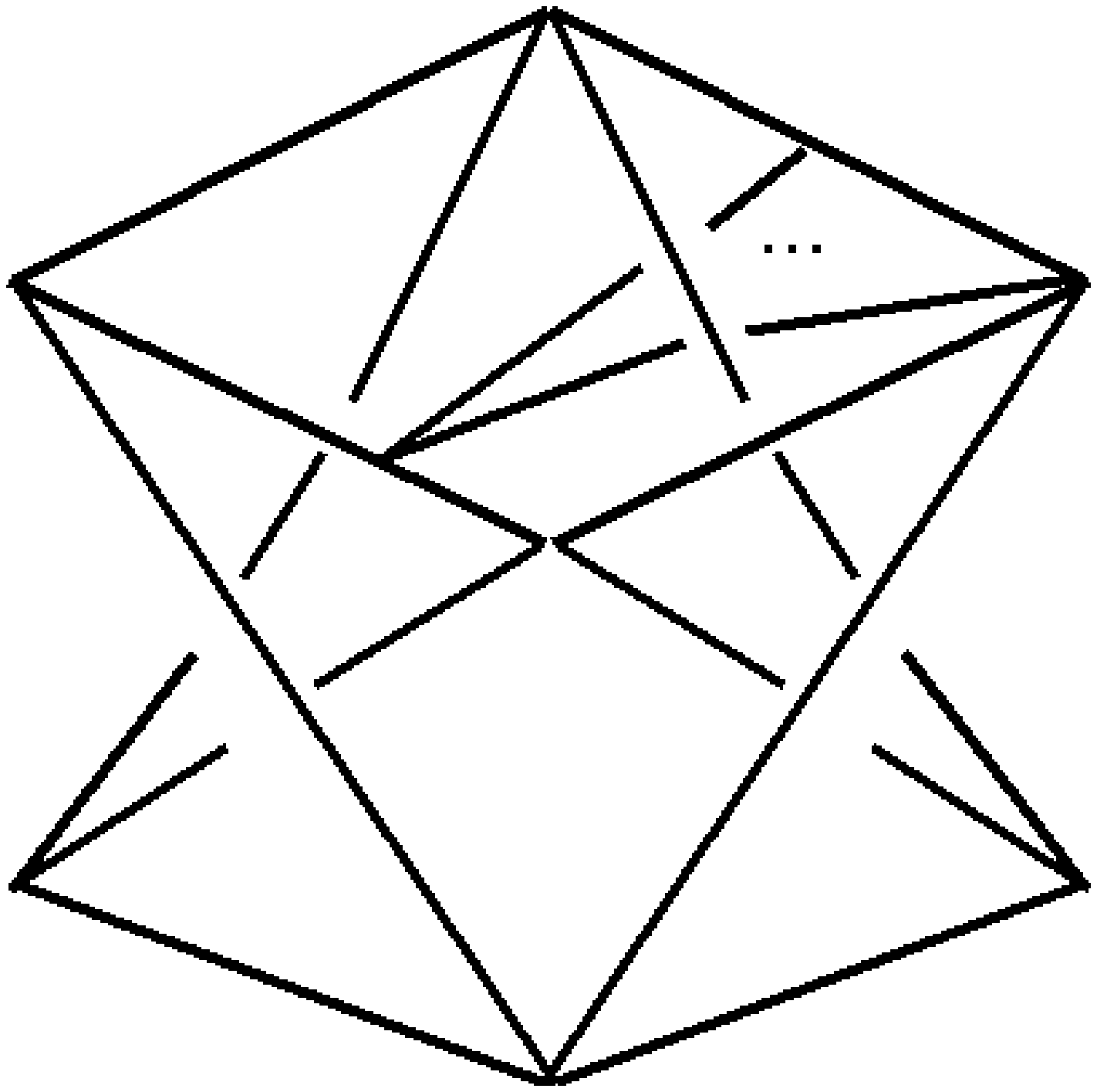}}\;}
\def\bip{\raisebox{-2cm}{\epsfysize=4cm\epsfbox{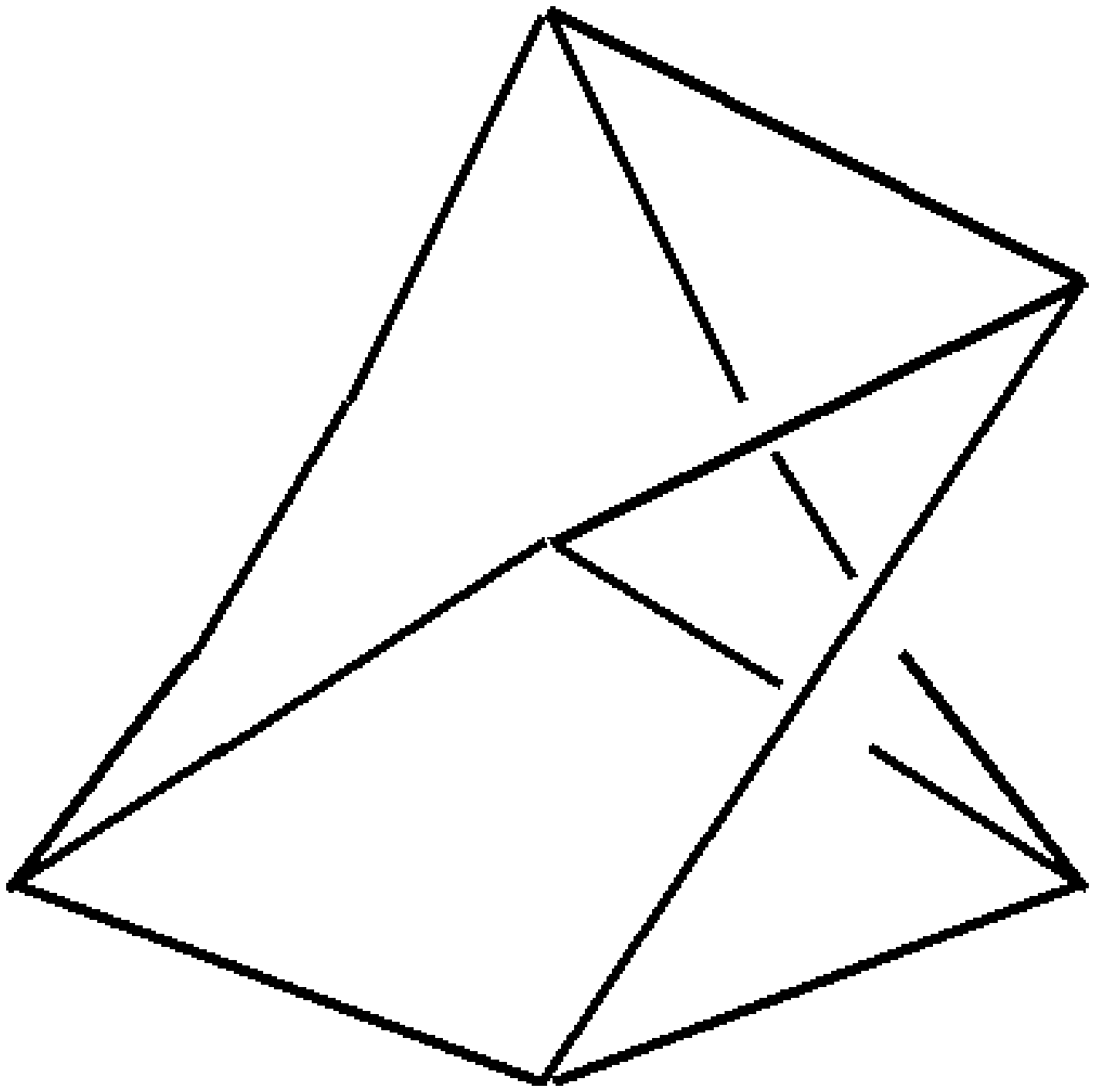}}\;}

\theoremstyle{plain}
\newtheorem{thm}{Theorem}
\newtheorem{lem}[thm]{Lemma}
\newtheorem{cor}[thm]{Corollary}
\newtheorem{prop}[thm]{Proposition}
\newtheorem{claim}[thm]{Claim}

\theoremstyle{definition}
\newtheorem{defn}[thm]{Definition}
\newtheorem{rmk}[thm]{Remark}

\newtheorem{ex}[thm]{Example}

\newtheorem{properties}[thm]{Properties}

\numberwithin{thm}{section}
\numberwithin{equation}{section}

\newcommand{\eq}[2]{\begin{equation}\label{#1}#2 \end{equation}}
\newcommand{\ml}[2]{\begin{multline}\label{#1}#2 \end{multline}}
\newcommand{\ga}[2]{\begin{gather}\label{#1}#2 \end{gather}}

\newcommand{\surj}{\twoheadrightarrow}
\newcommand{\inj}{\hookrightarrow}

\newcommand{\rank}{{\rm rank}}

\newcommand{\im}{{\rm im}}
\newcommand{\Spec}{{\rm Spec \,}}

\newcommand{\sF}{{\mathcal F}}

\newcommand{\sO}{{\mathcal O}}

\newcommand{\sQ}{{\mathcal Q}}

\newcommand{\sV}{{\mathcal V}}

\newcommand{\A}{{\mathbb A}}

\newcommand{\C}{{\mathbb C}}

\newcommand{\F}{{\mathbb F}}
\newcommand{\G}{{\mathbb G}}

\renewcommand{\P}{{\mathbb P}}
\newcommand{\Q}{{\mathbb Q}}
\newcommand{\R}{{\mathbb R}}

\newcommand{\Z}{{\mathbb Z}}

\begin{document}

\title[Graph Polynomials]{
On Motives Associated to Graph Polynomials}
\author{Spencer Bloch}
\address{Dept. of Mathematics,
University of Chicago,
Chicago, IL 60637,
USA}
\email{bloch@math.uchicago.edu}

\author{H\'el\`ene Esnault}
\address{Mathematik,
Universit\"at Duisburg-Essen, FB6, Mathematik, 45117 Essen, Germany}
\email{esnault@uni-essen.de}

\author{Dirk Kreimer}
\address{IHES,
91440 Bures sur Yvette, France and Boston U., Boston MA 02215.}
\email{kreimer@ihes.fr}
\date{Sept.\ 28b, 2005}

\begin{abstract}
The appearance of multiple zeta values in anomalous dimensions and
$\beta$-functions of renormalizable quantum field theories has
given evidence towards a motivic interpretation of these
renormalization group functions. In this paper we start to hunt
the motive, restricting our attention to a subclass of graphs in
four dimensional scalar field theory which give scheme independent
contributions to the above functions.
\end{abstract}
\subjclass{}
\maketitle
\begin{quote}

\end{quote}
\setcounter{section}{-1}
\section{Introduction}\label{secint}
Calculations of Feynman integrals arising in perturbative quantum
field theory \cite{Pisa,BK} reveal interesting patterns of zeta
and multiple zeta values. Clearly, these are motivic in origin,
arising from the existence of Tate mixed Hodge structures with
periods given by Feynman integrals. We are far from a detailed
understanding of this phenomenon. An analysis of the problem leads
via the technique of Feynman parameters \cite{ItzZub} to the study
of motives associated to graph polynomials. By the seminal work of
Belkale and Brosnan \cite{BB}, these motives are known to be quite
general, so the question becomes under what conditions on the
graph does one find mixed Tate Hodge structures and multiple zeta
values.

The purpose of this paper is to give an expository account of some
general mathematical aspects of these ``Feynman motives'' and to work out
in detail the special case of wheel and spoke graphs. We consider only
scalar field theory, and we focus on {\it primitively divergent} graphs.
(A connected graph $\Gamma$ is primitively divergent if $\# Edge(\Gamma) =
2h_1(\Gamma)$ where $h_1$ is the Betti number of the graph; and if
further for any connected proper subgraph the number of edges is strictly
greater than twice the first Betti number.) From a motivic point of view, these play
the role of ``Calabi-Yau'' objects in the sense that they have unique
periods. Physically, the corresponding periods are renormalization scheme
independent. 

Graph polynomials are introduced in sections \ref{sec1a} and \ref{sec1}
as special cases of discriminant polynomials associated to
configurations. They are homogeneous polynomials written in a preferred
coordinate system with variables corresponding to edges of the graph. The
corresponding hypersurfaces in projective space are {\it graph
hypersurfaces}. Section \ref{sec3} studies coordinate linear spaces
contained in the graph hypersurface. The normal cones to these linear
spaces are linked to graph polynomials of sub and quotient graphs.
Motivically, the chain of integration for our period meets the graph
hypersurface along these linear spaces, so the combinatorics of their
blowups is important. (It is curious that arithmetically interesting
periods seem to arise frequently (cf. multiple zeta values \cite{GM} or the
study of periods associated to Mahler measure in the non-expansive case
\cite{Den}) in situations where the polar locus of the integrand meets the
chain of integration in combinatorially interesting ways.) 

Section \ref{sect_global} is not used in the sequel. It exhibits a
natural resolution of singularities $\P(N) \to X$ for a graph
hypersurface $X$. $\P(N)$ is a projective bundle over projective space,
and the fibres $\P(N)/X$ are projective spaces.

Section \ref{sec4} introduces Feynman quadrics. The period of interest is
interpreted as an integral \eqref{4.3} over $\P^{2r-1}(\R)$. The
integrand has simple poles along $r$ distinct quadrics. When these
quadrics are associated to a graph $\Gamma$, the period is shown to be
convergent precisely when $\Gamma$ is primitively divergent as above. 

Section \ref{sec_schw} reinterprets the above period as a relative period
\eqref{5.10} associated to the graph hypersurface. This is the Schwinger
trick \cite{ItzZub}.

Section \ref{sec7} presents the graph motive in detail. Let $X\subset
\P^{2n-1}$ be the graph hypersurface associated to a primitive divergent
graph. Let $\Delta \subset \P^{2n-1}$ be the coordinate simplex (union of
$2n$ coordinate hyperplanes). An explicit sequence of blowups in
$\P^{2n-1}$ of linear spaces is described. Write $P \to \P^{2n-1}$ for the
resulting variety. Let $f:Y \subset P$ be the strict transform of $X$, and
let $B := f^{-1}(\Delta)$ be the total inverse image. Then the motive is
\eq{0.1}{H^{2n-1}(P\setminus Y, B\setminus B\cap Y)
} 

Section \ref{sect_motII} considers what can be said directly about the
motive of a graph hypersurface $X$ using elementary projection
techniques. The main tool is a theorem of C. L. Dodgson about determinants,
published in 1866. 

Section \ref{sec9} describes what the theory of motivic cohomology suggests
about graph motives in cases \cite{BK} where the period is related to a
zeta value. 

Section \ref{sec10} considers the Schwinger trick from a geometric point of
view. The main result is that in middle degree, the primitive cohomology of
the graph hypersurface is supported on the singular set. 

Sections \ref{sect_ws} and \ref{sec12} deal with wheel and spoke graphs. Write $X_n
\subset \P^{2n-1}$ for the hypersurface associated to the graph which is a
wheel with $n$ spokes. The main results are
\ga{}{H^{2n-1}_c(\P^{2n-1}\setminus X_n) \cong \Q(-2) \\
H^{2n-1}(\P^{2n-1}\setminus X_n) \cong \Q(-2n+3 ).
}
Further, the de Rham cohomology $H^{2n-1}_{DR}(\P^{2n-1}\setminus X_n)$ in this case
is generated by the integrand of our graph period \eqref{6.1}. Note that
nonvanishing of the graph period, which is clear by considerations of
positivity, only implies that the integrand gives a nonzero cohomology
class in $H^{2n-1}_{DR}(P\setminus Y, B\setminus B\cap Y)$. It does not a priori imply
nonvanishing in $H^{2n-1}_{DR}(P^{2n-1}\setminus X_n)$.

Finally, section \ref{sec13} discusses various issues which remain to be understood, including the
question of when the motive \eqref{0.1} admits a {\it framing}, the curious role of triangles in graphs
whose period is known to be related to a $\zeta$ value, and the possibility of constructing a Hopf
algebra $H$ of graphs such that assigning to a primitive divergent graph its motive would give rise to a
Hopf algebra map from $H$ to the Hopf algebra $MZV$ of mixed zeta values.

From a physics viewpoint, our approach starts with a linear algebra analysis of the
configurations given by a graph and its relations imposed by the
edges on the vertices, illuminating the structure of the graph
polynomial. An all important notion then is the one of a subgraph,
and the clarification of the correspondence between linear
subvarieties and subgraphs is our next achievement.

We then introduce the Feynman integral assigned to a Feynman
graphs based on the usual quadrics provided by the scalar
propagators of free field theory. The map from that Feynman
integral to an integration over the inverse square of the graph
polynomial proceeds via the Schwinger trick \cite{ItzZub}, which
we discuss in detail.

We next discuss the motive using relating chains of coordinate
linear subspaces of the graph hypersurfaces with chains of
subgraphs. This allows for a rather systematic stratification of
the graph hypersurface which can be carried through for the wheel
graphs, but fails in general. We give an example of such a
failure. The wheels are then subjected to a formidable computation
of their middle dimensional cohomology, a feast which we are at
the time of writing unable to repeat for even the next most simple
class of graphs, the zig-zag graphs of \cite{Pisa}, which, at each
loop order, evaluate indeed to a rational multiple of the wheel at
the same loop order. After collecting our results for the de Rham
class in the wheels case, we finish the paper with some outlook
how to improve the situation.
\\ \ \\
{\it Acknowledgement:} The second named author thanks Pierre Deligne for important discussions.

\section{Polynomials associated with Configurations}\label{sec1a}

Let $K$ be a field and let $E$ be a finite set. Write $K[E]$ for the
$K$-vector space spanned by $E$. A {\it configuration} is simply a linear
subspace  $i_V: V \inj K[E]$.  The space $K[E]$ is self-dual
in an evident way, so for $e
\in E$ we may consider the functional $e^\vee\circ i_V : V \to K$. Fix a
basis
$v_1,\dotsc,v_d$ for $V$, and let $M_e$ be the $d\times d$ symmetric
matrix associated to the rank $1$ quadratic form $(e^\vee\circ i_V)^2$
on $V$. Define a polynomial
\eq{1.1}{\Psi_V(A) = \det(\sum_{e\in E} A_eM_e).
}
$\Psi_V$ is homogeneous of degree $d$. Note that changing the basis of $V$
only changes $\Psi_V$ by a unit in $K^\times$.
\begin{rmk} Write $\iota_V: \P(V) \inj \P^{\#E-1}$ for the evident
embedding on projective spaces of lines. View the quadratic forms
$(e^\vee\circ i_V)^2$ as sections in $\Gamma(\P(V), \sO(2))$. Then
$\iota_V$ is defined by the possibly incomplete linear series spanned
by these sections, and $\Psi_V$ is naturally interpreted as defining the
{\it dual hypersurface} in $\P^{\#E-1,\vee}$ of sections of this linear
system which define singular hypersurfaces in $\P(V)$, cf. section
\ref{sect_global}.
\end{rmk}

\begin{lem}\label{lem1.2} Each $A_e$ appears with degree $\le 1$ in
$\Psi_V$.
\end{lem}
\begin{proof} The matrix $M_e$ has rank $\le 1$. If $M_e = 0$ then of
course $A_e$ doesn't appear and there is nothing to prove. If rank $M_e$
is $1$, then multiplying on the left and right by invertible matrices
(which only changes $\Psi_V$ by an element in $K^\times$) we may assume
$M_e$ is the matrix with $1$ in position $(1,1)$ and zeroes elsewhere. In
this case
\eq{}{\Psi_V = \det\begin{pmatrix}A_e+m_{ee} & \hdots \\
\vdots & \vdots \end{pmatrix}
}
where $A_e$ appears only in entry $(1,1)$. The assertion of the lemma
follows by expanding the determinant along the first row.
\end{proof}

As a consequence, we can write
\eq{}{\Psi_V(A) = \sum_{\{e_1,\dotsc,e_d\}}c_{e_1,\dotsc,e_d}
A_{e_1}A_{e_2}\cdots A_{e_d} }
\begin{lem}\label{lem1.3} With notation as above, write
$M_{e_1,\dotsc,e_d}$ for the matrix (with respect to the chosen basis of
$V$) of the composition
\eq{}{V \to K[E] \xrightarrow{e'\mapsto 0,\ e' \neq e_i} Ke_1\oplus \ldots
\oplus Ke_d.
}
Then $c_{e_1,\dotsc,e_d} = \det M_{e_1,\dotsc,e_d}^2$.
\end{lem}
\begin{proof}As a consequence of Lemma \ref{lem1.2},
$c_{e_1,\dotsc,e_d}$ is obtained from $\Psi_V$ by setting $A_{e_i}
= 1,\ 1\le i \le d$ and $A_{e'} = 0$ otherwise, i.e.
$c_{e_1,\dotsc,e_d} = \det(\sum_i M_{e_i})$. With respect to the chosen
basis of
$V$ we may write $e^\vee \circ i_V = \sum a_{e,i}v_i^\vee : V \to K$. Then
$M_e = (a_{e,i}a_{e,j})_{ij}$ so
\eq{}{M_{e_1,\dotsc,e_d} = (a_{e,i});\quad \sum_e M_e =
(a_{e,i})(a_{j,e})^t = M_{e_1,\dotsc,e_d}M_{e_1,\dotsc,e_d}^t.
}
\end{proof}

\begin{cor} The coefficients of $\Psi_V$ are the squares of the Pl\"ucker
coordinates for $K[E] \surj W$. More precisely, the coefficient of
$\prod_{e\not\in T} A_e$ is Pl\"ucker$_T(W)^2$.
\end{cor}

\begin{rmk}Let $G$ denote the Grassmann of all $V_d \subset K[E]$. $G$
carries a line bundle $\sO_G(1)\cong \det(\sV)^\vee$, where $\sV \subset
K[E]\otimes_K \sO_G$ is the universal subbundle. Sections of $\sO_G(1)$
arise from the dual map $\bigwedge^d K[E] \cong \Gamma(G,\det
\sV^\vee)$.   Lemma \ref{lem1.3} can be interpreted universally as
defining a section
\eq{1.6}{\Psi \in \Gamma(G\times \P(K[E]), \sO_G(2)\boxtimes \sO_\P(1)).
}
\end{rmk}

Define $W = K[E]/V$ to be the cokernel of $i_V$. Dualizing yields an
exact sequence
\eq{}{0 \to W^\vee \xrightarrow{i_{W^\vee}} K[E] \to V^\vee \to 0
}
and hence a polynomial $\Psi_{W^\vee}(A)$ which is homogeneous of degree
$\# E - d$.
\begin{prop} We have the functional equation
\eq{}{\Psi_V(A) = c\cdot (\prod_{e\in E} A_e)\Psi_{W^\vee}(A^{-1});\quad c
\in K^\times }
\end{prop}
\begin{proof} For $T \subset E$ with $\#T = \#E -d$, consider the diagram
\eq{}{\begin{CD}@. @. 0 \\
@. @. @VVV \\
@. @. K[T] @>\beta_T>> W \\
@. @. @VVV @| \\
0 @>>> V @>>> K[E] @>>> W @>>> 0 \\
@. @| @VVV \\
@. V @>\alpha_{E-T}>> K[E-T] \\
@. @. @VVV \\
@. @. 0
\end{CD}
}
Fix bases for $V$ and $W$ so the isomorphism $\det K[E] \cong \det
V\otimes \det W$ (canonical up to $\pm 1$) is given by $c \in K^\times$.
Then $c = \det \alpha_{E\setminus T}\det\beta_{T}^{-1}$. By the above, the
coefficient in $\Psi_V$ of $\prod_{e\not\in T} A_e$ is $\det
\alpha_{E\setminus T}^2$ while the coefficient of $\prod_{e\in T} A_e$ in
$\Psi_{W^\vee}$ is $(\det\beta_{T}^t)^2$. The proposition follows
immediately.
\end{proof}

\begin{rmk}Despite the simple relation between $\Psi_V$ and
$\Psi_{W^\vee}$ it is useful to have both. When we apply this machinery
in the case of graphs, $\Psi_{W^\vee}$ admits a much more concrete
description. On the other hand, $\Psi_V$ is more closely related to the
Feynman integrals and periods of motives.
\end{rmk}
\begin{rmk}Let $K[E] \surj W$ be as above, and suppose $W$ is given with a basis. Then the
matrix $\sum_e A_eM_e$ associated to $i_{W^\vee}: W^\vee \inj K[E]$ is canonical as well.
In fact, a situation which arises in the study of graph polynomials is an exact sequence
$K[E] \to W \to K \to 0$. In this case, the matrix $\sum A_e M_e$ has zero determinant.
Define $W^0 := {\rm Image}(K[E] \to W)$. It is easy to check that the graph polynomial for
$i_{W^{0\vee}}: W^{0\vee} \inj K[E]$ is obtained from $\sum A_e M_e$ by removing the first row
and column and taking the determinant.
\end{rmk}

\section{Graph Polynomials}\label{sec1}

A finite graph $\Gamma$ is given with edges $E$ and vertices $V$. We orient the
edges. Thus each vertex of $\Gamma$ has entering edges and exiting edges. For a given vertex $v$ and a given edge $e$, we define ${\rm sign}(v,e)$ to be
$-1 $ if $e$ enters $v$ and $+1$ if $e$ exists $v$. We  associate to $\Gamma$ a configuration (defined over $\Z$) via the homology sequence
\eq{2.1}{0 \to H_1(\Gamma,\Z) \to \Z[E] \xrightarrow{\partial} \Z[V] \to H_0(\Gamma,\Z) \to
0, }
where the bounday map is $\Z$-linear and defined by $\partial (e)=\sum_{v\in V}
{\rm sign}(v,e)\cdot v$. Then $\partial$ depends on the chosen orientation but
$H_i(\Gamma, \Z)$ do not.

When $\Gamma$ is connected, we write $\Z[V]^0 := \ker(\Z[V] \xrightarrow{\deg} \Z)$. We
define the {\it graph polynomial} of $\Gamma$
\eq{2.2}{\Psi_\Gamma := \Psi_{H_1(\Gamma, \Z)}.
}

Recall a {\it tree} is a connected and simply connected graph. A tree $T \subset \Gamma$ is
said to be a {\it spanning tree} for the connected graph $\Gamma$ if every
vertex of $\Gamma$ lies in $T$. (If $\Gamma$ is not connected, we can
extend the notion of spanning tree $T \subset \Gamma$ by simply requiring
that $T\cap \Gamma_i$ be a spanning tree in $\Gamma_i$ for each connected
component $\Gamma_i \subset \Gamma$.)

\begin{lem} Let $T$ be a subgraph of a connected graph $\Gamma$. Let $E =
E_\Gamma$ be the set of edges of $\Gamma$ and let $E_T \subset E$ be the
edges of $T$. Then $T$ is a spanning tree if and only if one has an exact
homology diagram as indicated:
\eq{2.3}{\begin{CD} @. @. 0 @. 0 \\
@. @. @VVV @VVV \\
@. @. \Z[E_T] @>\beta >\cong > \Z[V]^0 \\
@. @. @VVV @VVV \\
0 @>>> H_1(\Gamma) @>>> \Z[E] @>\partial >> \Z[V] @>>> \Z @>>> 0 \\
@. @| @VVV @VVV @| \\
0 @>>> H_1(\Gamma) @>\alpha >\cong > \Z[E\setminus E_T] @>> 0 > \Z @>>\cong > \Z @>>> 0 \\
@. @. @VVV @VVV \\
@. @. 0 @. 0
\end{CD}
}
\end{lem}
\begin{proof} Straightforward. \end{proof}

\begin{prop}\label{prop2.2} With notation as above, we have
\eq{2.4}{\Psi_\Gamma(A) = \sum_{T \text{span tr.}} \prod_{e\not\in T} A_e.
}
\end{prop}
\begin{proof} Fix a basis $h_j$ for $H_1(\Gamma)$. Then
\eq{2.5}{\Psi_\Gamma(A) = \det\Big(\sum_e A_e e^\vee(h_j)e^\vee(h_k)\Big)
}
Let $B \subset E$ have $b$ elements, and let $E' = E\setminus B$. The coefficient of the monomial
$\prod_{e\in B} A_e$
in $\Psi_\Gamma(A)$ is computed by setting $A_{e'} = 0$ for $e' \in E'$.
The coefficient is non-zero iff the
determinant
\eqref{1.1} is non-zero under this specialization,
and this is true iff we get a diagram as in \eqref{2.3}, i.e. iff $E' = E_T$ for a spanning
tree $T$.  The coefficient of this monomial is $1=\det(\alpha\alpha^t)$
where $\alpha$ is as in the bottom row of \eqref{2.3}. \end{proof}

\begin{rmk}\label{rmk2.3} If $\Gamma = \coprod \Gamma_i$ with $\Gamma_i$ connected, then
\eq{2.6}{\Psi_\Gamma = \prod_i \Psi_{\Gamma_i}
}
as both the free abelian group on edges and $H_1$ are additive in
$i$. If we define spanning ``trees'' in disconnected graphs as suggested above,
Proposition \ref{prop2.2} carries over to the disconnected case.
\end{rmk}
\begin{cor}\label{cor2.4} The coefficients of $\Psi_\Gamma$ are all either $0$ or $+1$.
\end{cor}
\begin{defn}\label{defn2.5}The graph hypersurface $X_\Gamma \subset
\P^{\#(E_\Gamma)-1}$ is the hypersurface cut out by $\Psi_\Gamma = 0$.
\end{defn}
\begin{properties}\label{p2.5}We list a certain  evident
properties of $\Psi_\Gamma$. \newline\noindent
1. $\Psi_\Gamma$ is a sum of monomials with coefficient $+1$.\newline\noindent
2. No variable $A_i$ appears with degree $>1$ in any monomial.\newline\noindent
3. Let $\Gamma_1$ and $\Gamma_2$ be graphs, and fix vertices $v_i \in \Gamma_i$. Define
$\Gamma := \coprod \Gamma_i/\{v_1\sim v_2\}$. Thus, $E_\Gamma = E_{\Gamma_1}\amalg
E_{\Gamma_1}$ and $H_1(\Gamma) = H_1(\Gamma_1)\oplus H_1(\Gamma_2)$. Writing $A^{(i)}$ for
the variables associated to edges of $\Gamma_i$, we see that $\Psi_\Gamma =
\Psi_{\Gamma_1}(A^{(1)})\Psi_{\Gamma_2}(A^{(2)})$. Geometrically, the graph hypersurface
$X_\Gamma : \Psi_\Gamma = 0$ is simply the join of the graph hypersurfaces $X_{\Gamma_i}$.
(Recall, if $P_i \subset \P^N$ are linear subsets of projective space such that $P_1\cap
P_2 = \emptyset$ and $\dim P_1 + \dim P_2 = N-1$, and is $X_i \subset P_i$ are closed
subvarieties, then the join $X_1*X_2$ is simply the union of all lines joining points of
$X_1$ to points of $X_2$.) In particular, if $\Gamma_2$ is a tree, so $\Psi_{\Gamma_2}
= 0$, then $X_\Gamma$ is a cone over $X_{\Gamma_2}$. \newline\noindent
4. Defining $\Psi_\Gamma$ via spanning trees \eqref{2.4} can lead to confusion in
degenerate cases. For example, if $\Gamma$ has only a single vertex ({\it tadpole
graph}) and $n$ edges, then $H_1(\Gamma) \cong \Z[E_\Gamma] \cong \Z^n$. Thus
$\Psi_\Gamma = \prod_1^n A_i$, but there are no spanning trees.
\end{properties}

\section{Linear Subvarieties of Graph Hypersurfaces}\label{sec3}

Let $\Gamma$ be a graph with $n = \# E_\Gamma$ edges. For convenience we take $\Gamma$ to
be connected. It will be convenient to use the notation $h_1(\Gamma) := \rank \
H_1(\Gamma)$. In talking about subgraphs of a given graph $\Gamma$, we will frequently not
distinguish between the subgraph and the collection of its edges. (In particular, we will
not permit isolated vertices.)

Recall we have associated to $\Gamma$ a hypersurface $X_\Gamma \subset
\P^{n-1}$. Our projective space has a distinguished set of homogeneous coordinates $A_e
\leftrightarrow e\in E_\Gamma$, so we get a dictionary
\ga{3.1}{\text{Subgraphs $G\subset \Gamma$} \leftrightarrow \text{coordinate linear
subspaces
$L
\subset \P^{n-1}$} \\
G \mapsto L(G):A_e=0, e \in G \notag \\
L: A_e=0, e\in S\subset E_\Gamma \mapsto G(L) = \bigcup_{e\in S} e \subset \Gamma. \notag
}

The Feynman period is the integral of a differential form on $\P^{n-1}$ with poles
along $X_\Gamma$ over a chain which meets $X_\Gamma$ along the non-negative real
loci of coordinate linear spaces contained in $X_\Gamma$. To give motivic meaning
to this integral, it will be necessary to blow up such linear spaces. The basic
combinatorial observation is

\begin{prop}\label{prop3.1} With notation as above, a coordinate linear space $L$ is
contained in
$X_\Gamma$ if and only if $h_1(G(L))>0$.
\end{prop}
\begin{proof}Suppose $L:A_e=0, e\in S$. Then $L \subset X_\Gamma$ if and only if every
monomial in $\Psi_\Gamma$ is divisible by $A_e$ for some $e\in S$. In other words, iff no
spanning tree of $\Gamma$ contains $S$. The assertion now follows from
\begin{lem}\label{lem3.2} Let $S \subset \Gamma$ be a (not necessarily connected) subgraph.
Then
$S$ is contained in some spanning tree for $\Gamma$ iff $h_1(S) = 0$.
\end{lem}
\begin{proof}[Proof of Lemma] Consider the diagram
\eq{}{\begin{CD} 0 @>>> H_1(S) @>>> \Z[E_S] @>c>> \Z[V_\Gamma]^0 \\
@. @VV i V @VV b V @| \\
0 @>>> H_1(\Gamma) @>a>> \Z[E_\Gamma] @>>> \Z[V_\Gamma]^0 @>>> 0.
\end{CD}
}
Note that the map $i$ is always injective. $S$ is itself a spanning tree iff $c$ is
surjective and $a$ and $b$ have disjoint images. If we simply assume disjoint images with
$c$ not surjective, we can find $e \in E_\Gamma$ such that $e\not\in \im(a)+\im(b)$. Then
$S' = S\cup \{e\}$ still satisfies $h_1(S')=0$. Continuing in this way, eventually $c$
must be surjective. Since the images of $a$ and $b$ remain disjoint, $c$ will be an
isomorphism, and the resulting subgraph of $\Gamma$ will be a spanning tree. \end{proof}

This completes the proof of the proposition. \end{proof}

Let $\Gamma$ be a connected graph as above, and let $G\subset \Gamma$ be a subgraph. It
will be convenient not to assume $G$ connected. In particular, $\Psi_G$ and $X_G$ will be
defined as in Remark \ref{rmk2.3}.  We define a modified quotient graph
\eq{3.3}{\Gamma \surj \Gamma//G
}
by identifying the connected components $G_i$ of $G$ to vertices $v_i \in \Gamma//G$
 (but not identifying $v_i\sim v_j$). If $G$ is connected, this is the standard
quotient in topology. One gets a diagram with exact rows and columns
\eq{3.4}{\begin{CD}@. 0 @. 0 @. 0\\
@. @VVV @VVV @VVV \\
0 @>>> H_1(G) @>>> \Z[E_G] @>>> \Z[V_G]^0 @>>> 0 \\
@. @VVV @VVV @VVV \\
0 @>>> H_1(\Gamma) @>>> \Z[E_\Gamma] @>>> \Z[V_\Gamma]^0 @>>> 0 \\
@. @VV \pi V @VVV @VVV \\
0 @>>> H_1(\Gamma//G) @>>> \Z[E_{\Gamma//G}] @>>> \Z[V_{\Gamma//G}]^0 @>>> 0 \\
@. @VVV @VVV @VVV \\
@. 0 @. 0 @. 0.
\end{CD}
}
Note with this modified quotient the map labeled $\pi$ is surjective.

Our objective now is to relate the graph hypersurfaces $X_\Gamma, X_G, X_{\Gamma//G}$. To
this end, we first consider the relation between spanning trees for the three graphs. If
$T\subset \Gamma$ is a spanning tree, then $h_1(T\cap G)=0$, but $T\cap G$ is not
necessarily connected. In particular it is not necessarily a spanning tree for $G$.

There is an evident lifting from subgraphs $V \subset \Gamma//G$ to subgraphs $\widetilde
V \subset \Gamma$ such that $\widetilde V$ and $G$ have no common edges.

\begin{lem}\label{lem3.3} Let $U \subset G$ be a spanning tree (cf. Remark \ref{rmk2.3}).
Then the association
\eq{3.5}{V \mapsto T:= \widetilde V \amalg U
}
induces a $1$ to $1$ correspondence between spanning trees $V$ of $\Gamma//G$ and spanning
trees $T$ of $\Gamma$ such that $U \subset T$.
\end{lem}
\begin{proof}Let $T$ be a spanning tree for $\Gamma$ and assume $U \subset T$.
Necessarily, $G\cap T=U$. Indeed, $U \subset G\cap T$ and $h_1(G\cap T)=0$. Since $U$ is
already a spanning tree, it follows from \eqref{2.3} that $G\cap T$ cannot be strictly
larger than $U$.

By \eqref{3.4}, $\pi(T)\cong T//U \subset \Gamma//G$ is connected and
$h_1(\pi(T)) = 0$. It follows that $\pi(T)$ is a spanning tree for $\Gamma//G$. We have $T
= \widetilde{\pi(T)}\amalg U$, so the association $T \mapsto \pi(T)$ is injective.

Finally, if $V \subset \Gamma//G$ is a spanning tree, then since
$$V \cong (\widetilde V\amalg U)//U,$$
it follows from \eqref{3.4} that $h_1(\widetilde V\amalg U)=0$. One easily checks that
this subgraph is connected and contains all the vertices of $\Gamma$, so it is a spanning
tree.
\end{proof}

\begin{prop}\label{prop3.4} Let $\Gamma$ be a connected graph, and let $G \subset \Gamma$
be a subgraph. Assume $h_1(G) = 0$. Let $X_\Gamma \subset \P(E_\Gamma)$ be the graph
hypersurface, and let $L(G):A_e=0, e\in G$ be the linear subspace of $\P(E_\Gamma)$
corresponding to $G$. Then $L(G)$ is naturally identified with $\P(E_{\Gamma//G})$, and
under this identification,
$$X_{\Gamma//G} = X_\Gamma \cap L(G).
$$
\end{prop}
\begin{proof} In this case, Lemma \ref{lem3.3} implies that spanning trees for $\Gamma//G$
are in $1$ to $1$ correspondence with spanning trees for $\Gamma$ containing $G$. It follows
from Proposition \ref{prop2.2} that
$$\Psi_{\Gamma//G} = \Psi_\Gamma |_{A_e=0,e \in G}.
$$
\end{proof}
\begin{prop}\label{prop3.5} Let $G\subset \Gamma$ be a subgraph, and suppose $h_1(G)>0$.
Then
$L(G):A_e=0, e\in G$ is contained in $X_\Gamma$. Let $P \to \P(E_\Gamma)$ be the blowup of
$L(G)
\subset \P(E_\Gamma)$, and let $F \subset P$
 be the exceptional locus. Let $Y \subset P$ be the strict transform of $X_\Gamma$ in $P$.
Then we have canonical identifications
\ga{}{ F \cong \P(E_G) \times \P(E_{\Gamma//G}) \\
Y \cap F = \Big(X_G\times \P(E_{\Gamma//G})\Big) \cup \Big(\P(E_G) \times X_{\Gamma//G}
\Big).
}
\end{prop}
\begin{proof}Let $T \subset \Gamma$ be a spanning tree. We have $h_1(T\cap G) = 0$ so
$T\cap G$ is contained in a spanning tree for $G$ by Lemma \ref{lem3.2}. In
particular, $\#(T\cap G) \ge \# E_G - h_1(G)$, with equality if and only if $T\cap G$ is a
spanning tree for $G$.

The normal bundle for $L(G) \subset \P(E_\Gamma)$ is $\bigoplus_{e\in G} \sO(1)$, from
which it follows that $F \cong L(G) \times \P(E_G)$. Also, of course, $L(G) \cong
\P(E_\Gamma \setminus E_G) \cong \P(E_{\Gamma//G})$.

We have $L(G) \subset X_\Gamma$ by Proposition \ref{prop3.1}. The intersection $F\cap Y$
is the projectivized normal cone of this inclusion. Algebraically, we identify
\eq{}{K[A_e]_{e\in \Gamma//G} \otimes K[A_e]_{e\in G}
}
with the tensor of the homogeneous coordinate rings for $\P(E_{\Gamma//G})$ and $\P(E_G)$.
Our cone is the hypersurface in this product defined by the sum of terms in $\Psi_\Gamma =
\sum_{T \subset \Gamma}\prod_{e \not\in T}A_e$ of minimal degree in the normal variables
$A_e, e\in G$. These correspond to spanning trees $T$ with $\# G\cap T$ maximal. By the
above discussion, these are the $T$ such that $T\cap G$ is a spanning tree for $G$. It now
follows from Lemma \ref{lem3.3} that in fact the cone is defined by
\eq{}{\Psi_{\Gamma//G}(A_e)_{e\in \Gamma//G} \cdot \Psi_{G}(A_e)_{e\in G} \in K[A_e]_{e\in
\Gamma//G} \otimes K[A_e]_{e\in G}.
}
The proposition is now immediate. \end{proof}

\begin{rmk}The set $F \cap Y$ above can also be interpreted as the exceptional fibre for
the blowup of $L(G) \subset X_\Gamma$.
\end{rmk}
\begin{ex}Fix an edge $e_0 \in \Gamma$ and take $G = \Gamma \setminus e_0$. Then $L(G) =:p$ is a
single point. If $p \not\in X_\Gamma$, then $h_1(G)=0$ and Proposition \ref{prop3.4}
implies that $X_{\Gamma//G} = \emptyset$. If $p \in X_\Gamma$ then $F \cong \P(E_\Gamma \setminus
e_0)$ and the exceptional divisor for the blowup of $p \in X_\Gamma$ is $X_{\Gamma\setminus e_0}$.

Algebraically, this all amounts to the identity
\eq{3.10}{\Psi_\Gamma = A_{e_0}\Psi_{\Gamma \setminus e_0} + \Psi_{\Gamma/e_0}
}
where the two graph polynomials on the right do not involve $A_{e_0}$.
\end{ex}

\section{Global Geometry} \label{sect_global}

In this section, for a vector bundle $E$ over a variety $X$ we write $\P(E)$ for the
projective bundle of hyperplane sections, so $a_*\sO_{\P(E)}(1)=E$, with $a: \P(E)\to X$.
In particular,  a
surjection of vector bundles $E \surj F$ gives rise to a closed immersion $\P(F) \inj
\P(E)$.

Consider projective space $\P^r$ and its dual $(\P^r)^\vee$. One has the Euler
sequence
\eq{}{0 \to \sO_{\P^r} \xrightarrow{e} \sO_{\P^r}(1) \otimes \Gamma((\P^r)^\vee, \sO(1))
\to T_{\P^r}
\to 0,
}
where $T$ is the tangent bundle. Writing $T_0,\dotsc,T_r$ for a basis of $\Gamma(\P^r,
\sO(1))$ and $\frac{\partial}{\partial T_i} \in \Gamma((\P^r)^\vee, \sO(1))$ for the dual
basis, we have
\eq{}{e(1) = \sum T_i \otimes \frac{\partial}{\partial T_i} \in
\Gamma\Big(\P^r,\sO_{\P^r}(1)
\otimes \Gamma((\P^r)^\vee, \sO(1))\Big).
}
Geometrically, we can think of $e(1)$ as a homogeneous form of degree $(1,1)$ on  $\P^r
\times (\P^r)^\vee$ whose zeroes define $\P(T_{\P^r}) \inj \P^r \times (\P^r)^\vee$. The
fibre in $\P(T_{\P^r})$ over a point $\frac{\partial}{\partial T_i} = a_i$ in
$(\P^r)^\vee$ is the hyperplane cut out by $\sum a_iT_i$ in $\P^r$.

For $V \inj \P^r$ a closed subvariety, define $p_V$ to be the composition
$p_V: \P(T_{\P^r}|_V) \inj \P(T_{\P^r}) \to (\P^r)^\vee$, and the fibre over
$\frac{\partial}{\partial T_i} = a_i$ is $V \cap \{\sum a_iT_i=0\}$. Assuming $V$ smooth,
we have the normal bundle sequence
\eq{}{0 \to T_V \to T_{\P^r}|_V \to N_{V/\P^r} \to 0.
}
\begin{prop}Assume $V \inj \P^r$ is a smooth, closed subvariety. Consider the diagram
\eq{}{\begin{CD} \P(N_{V/\P^r}) @>\inj >> \P(T_{\P^r}|_V) \\
@VVV @VV {p_V} V \\
 (\P^r)^\vee @=  (\P^r)^\vee
\end{CD}}
We have \eq{}{\P(N_{V/\P^r}) \cap p_V^{-1}(a) = (V \cap \{\sum
a_iT_i=0\})_{{\rm sing}}, } the singular points of the
corresponding hypersurface section.
\end{prop}
\begin{proof}Let $x \in V \subset \P^r$ be a point. To avoid confusion we write $dT_i$ for
the dual basis to $\frac{\partial}{\partial T_i}$. To a sum $\sum a_i dT_i$ and a point
$x\in V$ we can associate a point of $\P(T_{\P^r}|_V)$. Suppose $x \in p_V^{-1}(a)$. Then
$x$ is singular in this fibre if and only if $\sum a_i dT_i$ kills $T_{V,x} \subset
T_{\P^r,x}$, and this is true if and only if $\sum a_i dT_i \in \P(N_{V/\P^r})$.
\end{proof}

Suppose now $V = \P^k$ and the embedding $\P^k \inj \P^r$ is defined by a sublinear system
in $\Gamma(\P^k, \sO(2))$ spanned by quadrics $q_0,\dotsc,q_k$. The fibres of the map $p:
T_{\P^r/\P^k} \to (\P^r)^\vee$ are the degree $2$ hypersurfaces $\{\sum a_iq_i = 0\}
\subset \P^k$. Note that the singular set in such a hypersurface is a
projective space of dimension $= k - \rank(\sum a_i M_i)$, where the $M_i$ are $(k+1)\times
(k+1)$ symmetric matrices associated to the quadrics $q_i$. We conclude
\begin{prop} \label{prop44.2} With notation as above, define
\eq{}{X = \{ a \in (\P^r)^\vee | \rank (\sum a_iM_i) < k+1\}.
}
Then writing $N = N_{\P^k/\P^r}$, the map $\P(N) \to X$ is a resolution of singularities of $X$. The
fibres of this map are projective spaces, with general fibre $\P^0 = {\rm point}$.
\end{prop}

\section{Quadrics} \label{sec4}

Let $K \subset \R$ be a real field. (For the application to Feynman quadrics, $K=\Q$.) We
will be interested in homogeneous quadrics
\eq{}{Q_i : q_i(Z_1,\ldots,Z_{2r}) = 0, 1\le i\le r
}
in $\P^{2r-1}$ with homogeneous coordinates $Z_1,\dotsc,Z_{2r}$.
The union $\cup_i^r Q_i$ of the quadrics has then degree $2r$.  It implies that
$\Gamma(\P^{2r-1}, \omega(\sum_1^r Q_i))=K[\eta]$ for a generator $\eta$ which, on the affine open $Z_{2r}\neq 0$ with  affine coordinates $z_i=\frac{Z_i}{Z_r}, i=1,\dots, (2r-1)$, is $\eta|_{Z_{2r-1} \neq 0}=
\frac{dz_1\wedge \ldots \wedge dz_{2r-1}}{\tilde{q}_1\cdots \tilde{q}_r},$
with $\tilde q_i = \frac{q_i}{Z_{2r}^2}$.
By (standard) abuse of notations, we write
\eq{4.2}{\eta=\frac{\Omega_{2r-1}}{q_1\cdots q_{r}};\quad \Omega_{2r-1} :=
\sum_{i=1}^{2r} (-1)^i Z_idZ_1\wedge \cdots \widehat{dZ_i}\cdots \wedge dZ_{2r}. }
The
transcendental quantity of interest is the {\it period}
\eq{4.3}{P(Q) := \int_{\P^{2r-1}(\R)}\eta=
\int_{z_1,\dotsc,z_{2r-1}=-\infty}^\infty \frac{dz_1\wedge\cdots\wedge
dz_{2r-1}}{\tilde q_1\cdots \tilde q_{r}}.
}
The integral is convergent and the period well defined e.g. when the quadrics are all
positive definite.

Suppose now $r=2n$ above, so we consider quadrics in $\P^{4n-1}$. Let $H \cong K^n$ be a
vector space of dimension $n$, and identify $\P^{4n-1}=\P(H^4)$. For $\ell: H \to K$ a
linear functional, $\ell^2$ gives a rank $1$ quadratic form on $H$. A {\it Feynman
quadric} is a rank $4$ positive semi-definite form on
$\P^{4n-1}$ of the form $q=q_\ell = (\ell^2,\ell^2,\ell^2,\ell^2)$ . We will be interested
in quadrics $Q_i$ of this form (for a fixed decomposition $K^{4n} = H^4$.) In other words,
we suppose given linear forms $\ell_i$ on
$H$, $1\le i\le 2n$, and we consider the corresponding period $P(Q)$ where $q_i =
(q_{\ell_i},q_{\ell_i},q_{\ell_i},q_{\ell_i})$.

For $\ell: H \to K$ a linear form, write $\lambda = \ker(\ell),\ \Lambda =
\P(\lambda,\lambda,\lambda,\lambda) \subset \P(H^4)=\P^{4n-1}$. The Feynman quadric
$q_\ell$ associated to
$\ell$ is then a cone over the codimension $4$ linear space $\Lambda$. For a suitable
choice of homogeneous coordinates $Z_1,\dotsc,Z_{4n}$ we have $q_\ell = Z_1^2+\ldots +
Z_4^2$.

Let $q_1,\dotsc,q_{2n}$ be Feynman quadrics, and let $\Lambda_i$ be the linear space
associated to $q_i$ as above. As $K$ is a real field,  $\P^{4n-1}(\R)$ meets $Q_i(\C)$ only on
$\Lambda_i(\R)$.
\begin{lem}\label{lem4.1} With notation as above, for $I = \{i_1,\dotsc,i_p\} \subset
\{1,\dotsc,2n\}$, write $r(I) = {\rm codim}_H( \lambda_{i_1}\cap \ldots \cap
\lambda_{i_p})$. The integral
\eqref{4.2} converges if and only if $\sup_I \{ p(I) - 2r(I)\} < 0$. Here the sup is taken
over all $I\subset \{1,\dotsc,2n\}$ and $p(I) = \# I$.
\end{lem}
\begin{proof}Suppose $\lambda_1\cap\ldots\cap\lambda_p$ has codimension $r$, with $2r\le
p$. We can choose local coordinates $x_j$ so that $\bigcap_{i=1}^p \Lambda_i:x_1=\ldots
=x_{4r}=0$, and then make the blowup $y_j = \frac{x_j}{x_{4r}}, 1\le j\le 4r-1$, $y_j=x_j,\ j\ge
4r$. Then
\eq{}{\frac{d^{4n-1}x}{q_1(x)\cdots q_{2n}(x)} =
\frac{x_{4r}^{4r-1}d^{4n-1}y}{x_{4r}^{2p}\tilde q_1(y)\cdots \tilde q_{2n}(y)}
}
for suitable $\tilde q_i(y)$ which are regular in the $y$-coordinates. Since $|\prod
\tilde q_i^{-1}|\ge C>0$, it follows that the integral over a neighborhood of $0 \in
\R^{4n-1}$ diverges if $(4r-2p)\le 0$.

Suppose conversely that $\sup_I \{p(I)-2r(I)\}<0$. Note if $n=1$, the quadrics are smooth
and positive definite so the integrand has no pole along the integration chain and
convergence is automatic. Assume $n>1$. The above argument shows that blowing up an
intersection of the
$\Lambda_i$ does not introduce a pole in the integrand along the exceptional divisor.
Further, the strict transforms of the quadrics continue to have degree $\le 2$ in the
natural local coordinates and to be cones over the strict transforms of the $\Lambda_i$.
One knows  that after a finite number of such blowups, the strict
transforms of the $\Lambda_i$ will meet transversally (see \cite{ESV} for a minimal way to do it). All blowups and coordinates will be
defined over $K \subset \R$, and one is reduced to checking convergence for an integral of
the form
\eq{}{\int_U \frac{d^{4n-1}x}{(x_1^2+\ldots +x_4^2)\cdots (x_{4n-7}^2 + \ldots +
x_{4n-4}^2)}
}
with $U$ a neighborhood of $0 \in \R^{2n-1}$. The change of variables $x_i = ty_i, \ i\le
(4n-4)$ introduces a $t^{4n-5-2n+2}= t^{2n-3}$ factor. Since $n\ge 2$, convergence is clear.
\end{proof}

Let $\Gamma$ be a graph with $N$ edges and $n$ loops. Associated to $\Gamma$ we have the
configuration of $N$ hyperplanes in the $n$-dimensional vector space $H = H_1(\Gamma)$,
\eqref{2.1}. As above, we map consider the Feynman quadrics
$q_i=(\ell_i^2,\ell_i^2,\ell_i^2,\ell_i^2)$ on $\P^{4n-1}, 1\le i\le N$. The graph
$\Gamma$ is said to be {\it convergent} (resp. {\it logarithmically divergent}) if $N>2n$
(resp. $N=2n$). When $\Gamma$ is logarithmically divergent, the form
\eq{}{\omega_\Gamma := \frac{d^{4n-1}x}{q_1\cdots q_{2n}}
}
has poles only along $\bigcup Q_i$, and we define the period
\eq{4.7}{P(\Gamma):= \int_{\P^{4n-1}(\R)}\omega_\Gamma
}
as in \eqref{4.3}.
\begin{prop} \label{prop4.2} Let $\Gamma$ be a logarithmically divergent graph with $n$ loops and $2n$
edges. The period $P(\Gamma)$ converges if and only if every subgraph $G\subsetneq \Gamma$
is convergent, i.e. if and only if $\Gamma$ is primitive log divergent in the sense
discussed in section \ref{secint}.
\end{prop}
\begin{proof}Let $G \subset \Gamma$ be a subgraph with $m$ loops and $M$ edges, and
assume $M \le 2m$. Let $I \subset \{1,\dotsc,2n\}$ be the edges not in $G$. Note $H_1(G)
\subset H_1(\Gamma)$ has codimension $n-m$ and is defined by the $2n-M$ linear functionals
corresponding to edges in $I$. By Lemma \ref{lem4.1}, the fact that $2(n-m) \le 2n-M$
implies that the period integral $P(\Gamma)$ is divergent. Conversely, if the period
integral is divergent, there will exist an $I$ with $p(I)-2r(I)\ge 0$. Let $G \subset
\Gamma$ be the union of the edges not in $I$. Then $G$ has $2n-p(I)$ edges. Also $H_1(G)
\subset H_1(\Gamma)$ is defined by the vanishing of functionals associated to edges in
$I$, so $G$ has $n-r(I)$ loops. It follows that $G$ is not convergent.
\end{proof}

\section{The Schwinger Trick}\label{sec_schw}

Let $Q_i: q_i(Z_1,\dotsc,Z_{4n})=0, 1\le i\le 2n$ be quadrics in $\P^{4n-1}$. We assume
the period integral \eqref{4.3} converges. Let $M_i$ be the $4n\times 4n$ symmetric matrix
corresponding to $q_i$, and write
\eq{5.1}{\Phi(A_1,\dotsc,A_{2n}) := \det(A_1M_1+\ldots +A_{2n}M_{2n}).
}
The Schwinger trick relates the period integral $P(Q)$ \eqref{4.3} to an integral on
$\P^{2n-1}$
\eq{5.2}{\int_{\P^{4n-1}(\R)}\frac{\Omega_{4n-1}(Z)}{q_1\cdots q_{2n}} =
C\int_{\sigma^{2n-1}(\R)} \frac{\Omega_{2n-1}(A)}{\sqrt{\Phi}}.
}
Here $\sigma^{2n-1}(\R)\subset \P^{2n-1}(\R)$ is the locus of all points
$s=[s_1,\dotsc,s_{2n}]$ such that the projective coordinates $s_i \ge 0$. $C$ is an
elementary constant, and the $\Omega$'s are as in \eqref{4.2}. Note the homogeneity is such
that the integrands make sense.
\begin{lem}\label{lem5.1} With notation as above, define
\eq{5.3}{g(A) =
\int_{\P^{4n-1}(\R)}\frac{\Omega_{4n-1}}{(A_1q_1+\ldots+A_{2n}q_{2n})^{2n}}. }
Then
\eq{5.4}{g(A)\sqrt{\Phi} = c\pi^{-2n};\quad c \in \overline\Q^\times,\ [\Q(c):\Q] \le 2.
}
\end{lem}
If $\Phi = \Xi^2$ for a polynomial $\Xi \in \Q[A_1,\dotsc,A_{2n}]$, then $c\in \Q^\times$.
\begin{proof}By analytic continuation, we may suppose that $Q_a: \sum A_i q_i=0$ is smooth.
The integral is then the period associated to $H^{4n-1}(\P^{4n-1}\setminus Q_a)$. As generator for
the homology we may either take $\P^{4n-1}(\R)$ or the tube  $\tau \subset \P^{4n-1}\setminus Q_a$
lying over the difference of two rulings $\ell_1-\ell_2$ in the even dimensional smooth
quadric $Q_a$. (More precisely, let $S \subset N \xrightarrow{p} X$ be
the sphere bundle for some metric on the normal bundle $N$ of $X$,
where $X\subset \P^{2n-1}$ is defined by $\Phi=0$.  Take $\tau =
p^{-1}(\ell_1-\ell_2)$.) The two generators differ by a rational scale factor
$c$. Integrating over $\tau$ shows that $g(A)$ is defined up to a scale factor $\pm 1$ on
$\P^{2n-1}\setminus X$. The monodromy arises because the rulings $\ell_i$ on $Q_a$ can be
interchanged as $a$ winds around $X$. It follows easily that the left hand side in
\eqref{5.4} is homogeneous of degree $0$ and single-valued on $\P^{2n-1}\setminus X$. To study its
behavior near $X$ we restrict to a general line in $\P^{2n-1}$. In affine coordinates, we
can then assume the family of quadrics looks like $(\sum_1^{4n-1}x_i^2) - t=0$, where $t$ is
a parameter on the line. The integral then becomes
\eq{}{\int_\gamma \frac{dx_1\wedge\ldots\wedge dx_{4n-1}}{(\sum x_i^2 -t)^{2n}} = {\rm const}\cdot
t^{-\frac{1}{2}}
}
for a suitable cycle $\gamma$.
The change of variable $x_i=y_it^{\frac{1}{2}}$ gives the value ${\rm const}\cdot
 t^{-\frac{1}{2}}$ from which one
sees that $g(A)\sqrt{\Phi}$ is constant. Since $H^{4n-1}(\P^{4n-1}\setminus Q_a) \cong \Q(-2n)$ as
Hodge structure, $g(A) = c_0\pi^{-{2n}}$ for some $c_0 \in \Q^\times$, and the lemma
follows.
\end{proof}

With notation as above, define
\eq{}{f(A) := \int_{\P^{4n-1}(\R)} \frac{\Omega_{4n-1}(Z)}{(A_1q_1+\ldots
+A_{2n}q_{2n})q_2q_3\cdots q_{2n}}
}
Note that $f(A)$ is defined for $q_i$ positive definite and $A_j \ge 0$ but not all
$A_j=0$. We have
\eq{}{g(A) = \frac{-1}{(2n-1)!}\frac{\partial^{2n-1}}{\partial A_2\ldots \partial A_{2n}}
f(A).
}
Write $a_i = \frac{A_i}{A_1}, 2\le i\le 2n$, and define $F(a_2,\dotsc,a_{2n}):= A_1f(A)$. Note the
various partials $\partial^{i-1}/\partial a_2\ldots \partial a_iF(a)$
vanish as
$a_i \to +\infty$ with $a_j \ge 0, \forall j$. Also $\frac{\Omega_{2n-1}}{A_1^{2n}} =
-da_2\wedge\ldots \wedge da_{2n-1}$. Thus
\ml{5.8}{\int_{\sigma^{2n-1}(\R)} g(A)\Omega_{2n-1}(A) = -\int_{\sigma^{2n-1}(\R)}A_1^{2n}
g(A)da_2\wedge\ldots \wedge da_{2n-1} = \\
\frac{1}{(2n-1)!}\int_{a_2,\dotsc,a_{2n}=0}^{+\infty}\frac{\partial^{2n-1}}{\partial
a_2\ldots
\partial a_{2n}} F(a)da_2\wedge\ldots \wedge da_{2n-1} = \\
\frac{-1}{(2n-1)!}F(0,\dotsc,0) = \int_{\P^{4n-1}(\R)}\frac{\Omega_{4n-1}(Z)}{q_1q_2\cdots
q_{2n}} = P(Q).
}
This identity holds by analytic extension in the $q$'s where both integrals are defined.
Combining \eqref{5.8} with Lemma \ref{lem5.1} we conclude
\begin{prop}With notation as above, assuming the integral defining $P(Q)$ is convergent,
we have
\eq{}{P(Q) := \int_{\P^{4n-1}(\R)}\frac{\Omega_{4n-1}(Z)}{q_1q_2\cdots
q_{2n}} = \frac{c}{\pi^{2n}}\int_{\sigma^{2n-1}(\R)}\frac{\Omega_{2n-1}(A)}{\sqrt{\Phi}}.
}
\end{prop}
\begin{cor} \label{cor5.3} Let $\Gamma$ be a graph with $n$ loops and $2n$ edges. Assume every proper
subgraph of $\Gamma$ is convergent, and let $q_1,\dotsc,q_{2n}$ be the Feynman quadrics
associated to $\Gamma$ (cf. section \ref{sec4}). The symmetric matrices $M_i$ \eqref{5.1}
in this case are block diagonal
$$M_i = \begin{pmatrix}N_i & 0 & 0 & 0 \\ 0 &
N_i & 0 & 0 \\
0 & 0 & N_i & 0 \\
0 & 0 & 0 & N_i\end{pmatrix}$$
and $\Phi = \Psi_\Gamma^4$, where $\Psi_\Gamma = \det(A_1N_1+\ldots + A_{2n}M_{2n})$ is
the graph polynomial \eqref{2.2}. The Schwinger trick yields (cf. \eqref{4.7})
\eq{5.10}{P(\Gamma) := \int_{\P^{4n-1}(\R)}\frac{\Omega_{4n-1}(Z)}{q_1q_2\cdots
q_{2n}} = \frac{c}{\pi^{2n}}\int_{\sigma^{2n-1}(\R)}\frac{\Omega_{2n-1}(A)}{\Psi_\Gamma^2}.
}
for $c \in \Q^\times$.
\end{cor}

\section{The Motive}\label{sec7}

We assume as in section 5 that the ground field $K\subset \R$ is real.
Let $\Gamma$ be a graph with $n$ loops and $2n$ edges and assume every proper
subgraph of $\Gamma$ is convergent. Our objective in this section is to consider the
motive with period
\eq{6.1}{\int_{\sigma^{2n-1}(\R)}\frac{\Omega_{2n-1}(A)}{\Psi_\Gamma^2}.
}

We consider $\P^{2n-1}$ with fixed homogeneous coordinates $A_1,\dotsc,A_{2n}$ associated
with the edges of $\Gamma$. Linear spaces $L \subset \P^{2n-1}$ defined by vanishing of
subsets of the $A_i$ will be referred to as coordinate linear spaces. For such an $L$, we
write $L(\R^{\ge 0})$ for the subset of real points with non-negative coordinates.
\begin{lem}$X_\Gamma(\C) \cap \sigma^{2n-1}(\R) = \bigcup_{L \subset X_\Gamma}L(\R^{\ge 0})$,
where the union is taken over all coordinate linear spaces $L \subset X_\Gamma$.
\end{lem}
\begin{proof} We know by Corollary \ref{cor2.4} that $\Psi_\Gamma$ is a sum of monomials
with coefficients $+1$. The lemma is clear for the zero set of any polynomial with
coefficients $>0$.
\end{proof}
\begin{rmk}(i) The assertion of the lemma is true for any graph polynomial. We do not need
hypotheses about numbers of edges or loops. \newline\noindent
(ii) By Proposition \ref{prop3.1}, coordinate linear spaces $L \subset X_\Gamma$
correspond to subgraphs $G \subset \Gamma$ such that $h_1(G)>0$.
\end{rmk}

\begin{prop}\label{prop6.3} Let $\Gamma$ be as above. Define
\eq{}{\eta = \eta_\Gamma = \frac{\Omega_{2n-1}(A)}{\Psi_\Gamma^2}
}
as in \eqref{4.2}.
There exists a tower
\ga{6.3}{P = P_r \xrightarrow{\pi_{r,r-1}} P_{r-1} \xrightarrow{\pi_{r-1,r-2}} \ldots
\xrightarrow{\pi_{2,1}} P_1\xrightarrow{\pi_{1,0}} \P^{2n-1};\\ \pi =
\pi_{1,0}\circ\cdots\circ\pi_{r,r-1} \notag }
where $P_i$ is obtained from $P_{i-1}$ by blowing up the strict transform of a coordinate
linear space $L_i \subset X_\Gamma$ and such that \newline\noindent
(i) $\pi^*\eta_\Gamma$ has no poles along the exceptional divisors associated to the
blowups.\newline\noindent
(ii) Let $B \subset P$ be the total transform in $P$ of the union of coordinate
hyperplanes $\Delta^{2n-2} : A_1A_2\cdot A_{2n}=0$ in $\P^{2n-1}$. Then $B$ is a normal
crossings divisor in $P$. No face (= non-empty intersection of components) of $B$ is
contained in the strict transform $Y$ of $X_\Gamma$ in $P$. \newline\noindent
(iii) the strict transform of $\sigma^{2n-1}(\R)$ in $P$ does not meet $Y$.
\end{prop}
\begin{proof}Our algorithm to construct the blowups will be the following. Let $S$ denote
the set of coordinate linear spaces $L \subset \P^{2n-1}$ which are maximal, i.e. $L\in S,
L \subset L' \subset X_\Gamma \Rightarrow L=L'$. Define
\eq{}{\sF = \{ L \subset X_\Gamma \text{ coordinate linear space}\ |\ L=\bigcap L^{(i)},\
L^{(i)} \in S\}.
}
Let $\sF_{{\rm min}} \subset \sF$ be the set of minimal elements in $\sF$. Note that elements of
$\sF_{{\rm min}}$ are disjoint. Define $P_1 \xrightarrow{\pi_{1,0}} \P^{2n-1}$ to be the blowup
of elements of $\sF_{{\rm min}}$. Now define $\sF_1$ to be the collection of strict transforms
in $P_1$ of elements in $\sF\setminus \sF_{{\rm min}}$. Again elements in $\sF_{1,{\rm min}}$ are disjoint, and
we define $P_2$ by blowing up elements in $\sF_{1,{\rm min}}$. Then $\sF_2$ is the set of strict
transforms in $P_2$ of $\sF_1 \setminus  \sF_{1,{\rm min}}$, etc. This process clearly terminates.

Note that to pass from $P_i$ to $P_{i+1}$ we blow up strict transforms of coordinate
linear spaces $L$ contained in $X_\Gamma$. There will exist an open set $U \subset
\P^{2n-1}$ such that $P_{i} \times_{\P^{2n-1}} U \cong U$ and such that $L \cap U \neq
\emptyset$. It follows that to calculate the pole orders of $\pi^*\eta_\Gamma$ along
exceptional divisors arising in the course of our algorithm it suffices to consider the
simple blowup of a coordinate linear space $L\subset X_\Gamma$ on $\P^{2n-1}$. Suppose
$L:A_1=\ldots A_p=0$. By assumption, the subgraph $G = \{e_1,\dotsc,e_p\} \subset \Gamma$
is convergent, i.e. $p>2h_1(G)$. As in Proposition \ref{prop3.5}, if $I = (A_1,\dotsc,A_p)
\subset K[A_1,\dotsc,A_{2n}]$, then $\Psi_\Gamma \in I^{h_1(G)}- I^{h_1(G)+1}$ so the
denominator of $\eta_\Gamma$ contributes a pole of order $2h_1(G)$ along the exceptional
divisor. On the other hand, writing $a_i = \frac{A_i}{A_{2n}}$, a typical open in the blowup will
have coordinates $a_i' = \frac{a_i}{a_p}, i<p$ together with $a_p,\dotsc,a_{2n-1}$ and the
exceptional divisor will be defined by $a_p=0$. Thus
\ml{}{da_1\wedge\ldots\wedge da_{2n-1} = d(a_pa_1')\wedge \ldots \wedge
d(a_pa_{p-1}')\wedge da_p\wedge \ldots = \\
a_p^{p-1}da_1'\wedge\ldots\wedge da_{p-1}'\wedge da_p\ldots
}
Finally, $\pi^*\eta$ will vanish to order $p-1-2h_1(G)\ge 0$ on the exceptional divisor,
so the algorithm will imply (i). Here we observe that at least on the strata for which $p$ is even, $\pi^*\eta$ not only is regular along the exceptional divisor, but indeed really vanishes to order $\ge 1$.

Recall the dictionary \eqref{3.1} between subgraphs $G=G(L)\subset \Gamma$ and coordinate
linear spaces $L=L(G)$.

\begin{lem}\label{lem6.4} Let $\sF$ be as above, and let $\emptyset \neq L_1 \subsetneq
L_2 \subsetneq \ldots \subsetneq L_r$ be a chain of faces in $\sF$ which is saturated in
the sense that it cannot be made longer using elements of $\sF$. Let $G_r
\subsetneq G_{r-1}\ldots\subsetneq G_1\subsetneq G_0:=\Gamma$ be the chain
of subgraphs. Then $h_1(G_j) = r+1-j$. In particular, $n = h_1(\Gamma)
= r+1$. For $j\ge 1$ and any $e\in G_j \setminus  G_{j+1}$ we have $h_1(G_j\setminus e)=
h_1(G_j)-1 = h_1(G_{j+1})$.
\end{lem}
\begin{proof}[Proof of Lemma] Let $G \subset \Gamma$ be a (not necessarily connected)
subgraph. Consider the property
\eq{3}{\forall e \in G,\ h_1(G\setminus e)<
h_1(G).
}
I claim we can write $G = \bigcup G^{(i)}$ where the $G^{(i)}$
have the same minimality property and in addition $h_1(G^{(i)}) = 1$. We
argue by induction on $h=h_1(G)$. If $h=1$ we can just take $G$. If
$h>1$, then for every $e\in G$ we can find a $G_e \subset G$ such that
$e\in G_e$,
$h_1(G_e)=1$, and $G_e$ is minimal. Indeed, since $h_1(G\setminus e)<h_1(G)$, we
can find a connected subgraph $G' \subset G$ such that $e\in G'$,
$h_1(G')=1$, and $h_1(G'\setminus e)=0$. Now just remove $e' \neq e$ from $G'$
until the resulting subgraph is minimal.

Since $e\in G_e$ we have $G= \bigcup G_e$ as desired. Applying our
dictionary, $L(G) = \bigcap L(G_e)$. Note the $L(G_e) \subset X$ are
maximal. We conclude that $L(G) \in \sF$ for any $G\subset \Gamma$
satisfying \eqref{3}. Conversely, if $G = \bigcup G^{(i)}$ with
$L(G^{(i)})$ maximal in $X$, then every vertex in $G$ lies on at least
$2$ edges (because this holds for the $G^{(i)}$). If for some $e\in G$ we
had $h_1(G) = h_1(G\setminus e)$, we would then necessarily have that $G\setminus e$ was
disconnected. If $e \in G^{(1)} \subset G$, then since $G^{(1)}$ has no
external edges, it would follow that $G^{(1)}\setminus e$ was disconnected. This
would imply $h_1(G^{(1)}\setminus e) = h_1(G^{(1)})$, a contradiction.

We conclude that $L \in \sF$ iff $G(L)$ satisfies \eqref{3}. The lemma
now is purely graph-theoretic, concerning the existence of chains of subgraphs satisfying
\eqref{3}.  Basically the condition is that the $G_i$
have no external edges and are ``$1$-particle irreducible'' in the physicist's sense. (Note
of course that we cannot assume this for $G_0 = \Gamma$, which is given.) To construct
such a chain one simply takes $G_r \subset \Gamma$ minimal such that $h_1(G_r)=1$ and
$G_{r-i}$ minimal such that $G_{r-i+1} \subset G_{r-i}$ and $h_1(G_{r-i+1}) >
h_1(G_{r-i})$. Note the $G_j$ are not necessarily connected.
\end{proof}

We now prove (ii). Let $\pi:P \to \P^{2n-1}$ be constructed as above, using the
$\sF_{i,{\rm min}}$. $0$-faces of $B \subset P$ will be referred to as {\it vertices} (not to be
confused with vertices of the graph). It will suffice to show that no vertex lies in the
strict transform $Y$. Let $v \in P$ be a vertex. The question of whether $v \in Y$ is local
around $v$, so we may localize our tower \eqref{6.3}, replacing $P_i$ with
$\Spec(\sO_{P_i, v_i})$ where $v_i \in P_i$ is the image of $v$. In
particular, $\P^{2n-1}$ is replaced by $\Spec(\sO_{\P^{2n-1},v_0})$, where
$v_0\in \P^{2n-1}$ is the image of $v$. Note the image $v_i$ of $v$ in
$P_i$ is always a vertex.

We modify the tower by throwing out the steps for which $\Spec(\sO_{P_i,
v_i}) \to \Spec(\sO_{P_{i-1}, v_{i-1}})$ are isomorphisms. For
convenience, we don't change notation. All our $P_i$ are now local. Let
$E_1,\dotsc,E_r \subset P$ be the exceptional divisors, where $E_i$ comes
by pullback from $P_i$. Write $L_i := \pi(E_i) \subset P_0 :=
\Spec(\sO_{\P^{2n-1},v_0})$. We claim that $v_0 \in L_1$, and $L_1\subsetneq
L_2 \subsetneq \ldots \subsetneq L_r$ is precisely the sort of saturated
chain in $\sF$ considered in Lemma \ref{lem6.4} above. Indeed, at each
stage, $v$ maps to the exceptional divisor from the stage before. (If $v$
does not map to the exceptional divisor in $P_i$, then the local rings at
the image of $v$ in $P_i$ and $P_{i-1}$ are isomorphic, and this arrow is
dropped under localization.)

Our task now will be to compute $Y\cap\bigcap_{i=1}^r E_i$. We will do
this step by step. (We drop the assumption that our chain is
saturated.) Suppose first
$r=1$, i.e. there is only one blowup. Let $L_1\subset \P^{2n-1}$ be the
linear space being blown and suppose $L_1$ has codimension $p_1$. Then by Proposition
\ref{prop3.5} if we write
$G_1=G(L_1) \subset
\Gamma$ and
$\Gamma//G_1$ for the quotient identifying each connected component of $G$
to a point, we have $E_1 \cong L_1 \times \P^{p_1-1}$ and
\eq{6.7}{Y_1\cap E_1 = (X_{\Gamma//G_1} \times \P^{p_1-1}) \cup (L_1 \times
X_{G_1}).
}
Now suppose we have $L_1 \subset L_2$ and we want to compute $Y_2 \cap
E_1\cap E_2 \subset P_2$. (We write abusively $E_1$ for the pullback to
$P_2$ of $E_1$. $Y_i \subset P_i$ is the strict transform of $X$.). Locally at $v_0$ let
$L_i: a_1=\ldots = a_{p_i}=0$ with
$p_1>p_2$. Let $f$ be a local defining equation for $X$ near $v_0$ and
write
\eq{6.8}{f = \sum c_{I,J}(a_1,\dotsc,a_{p_2})^I(a_{p_2+1},\dotsc,a_{p_1})^J
}
with evident multi-index notation. Write $|I|, |J|$ for the
total degree of a multi-index. We are interested in points of  $P_1$
where the strict transform of $L_2$ meets $E_1$. Typical local
coordinates at such points look like
\eq{}{a_i':=a_i/a_{p_1},\ 1\le i < p_1,\ a_{p_1}'=a_{p_1},\ldots(\text{
coords. not involving the a's}).
}
To compute the intersection of the strict transform with the two
exceptional divisors on $P_2$, we let $\nu := \min(|I|+|J|)$ in
\eqref{6.8}, and write
\eq{8}{f_1 = \sum
(a_{p_1}')^{|I|+|J|-\nu}
c_{I,J}(a_1',\dotsc,a_{p_2}')^I(a_{p_2+1}',\dotsc,a_{{p_1}-1}')^J. }
This is the equation for $Y_1 \subset P_1$. We then take the image in the
cone for the second blowup by taking the sum only over those terms with
$|I|=|I|_{min}$ minimal:
\eq{6.11}{ \tilde f_1 = \sum_{\substack{I,J\\ |I|=|I|_{min}}}
(a_{p_1}')^{|I|+|J|-\nu}
c_{I,J}(a_1',\dotsc,a_{p_2}')^I(a_{p_2+1}',\dotsc,a_{{p_1}-1}')^J.
}
Notice that {\it a priori} $a_{p_1}'$ might divide $\tilde f_1$. We claim
in fact that it does not, i.e. that there exists $I,J$ such that $c_{I,J}
\neq 0$ and both $|I|$ and $|I|+|J|$ are minimum. To see this, note
\eq{}{|I|_{{\rm min}} = h_1(G_2);\quad \min(|I|+|J|) = h_1(G_1)
}
Assuming $L_1 \subset L_2$ is part of a saturated tower, we have as in
Lemma \ref{lem6.4} that $h_1(G_1) = h_1(G_2)+1$. If no nonzero term in $f$
has both $|I|$ and $|I|+|J|$ minimal, then every term with $|I|+|J|$
minimal must have $|I| = |I|_{{\rm min}} + 1$ and $|J|=0$. But this would mean
that the graph polynomial for $G_1$ would not involve the variables
$A_{p_2+1},\dotsc,A_{p_1}$. Since the $G_i$ have no external edges and
$h_1(G_i \setminus  e) < h_1(G_i)$, there are spanning trees ( disjoint unions of
spanning trees if $G_i$ is not connected) avoiding any given edge, so
this is a contradiction.

In general, if we have $L_1\subset\ldots\subset L_r$ saturated we write
\eq{6.13}{f = \sum_{I_1,\dotsc,I_r} c_{I_q,\dotsc,I_r}
(a_1,\dotsc,a_{p_r})^{I_r} (a_{p_r+1},\dotsc,a_{p_{r-1}})^{I_{r-1}}\cdots
(a_{p_2+1},\dotsc,a_{p_1})^{I_1}
}
We have
\ml{6.14}{\min(|I_r|) = \min(|I_{r-1}|+|I_r|) - 1 = \\
\ldots = \min(|I_r|+\cdots + |I_1|) -r+1
}
We claim there exist spanning trees $T$ for $G_1$ such that $T$ does not
contain any $G_i \setminus  G_{i+1}$. This will mean there exist $
c_{I_q,\dotsc,I_r} \neq 0$ such that $\sum_1^r |I_j|$ is minimum but
$|I_j| \neq 0$ for any $j$. By
\eqref{6.14}, this in turn implies for such a monomial that $\sum_{i=q}^r
|I_i|$ is minimal for all $q$. To show the existence of $T$, choose $e_i
\in G_i \setminus G_{i+1}$ for $1\le i\le r-1$ and $e_r \in G_r$. It suffices to
show that $h_0(G_1\setminus \{e_1,\ldots,e_r\}) = h_0(G_1)$. We have $h_0(G_1\setminus e_1) =
h_0(G_1)$ (since $h_1$ drops). Meyer Vietoris yields an exact sequence
\ml{6.15}{\ldots \to H_1(G_1\setminus e_1) \to
H_0(G_2\setminus \{e_2,\ldots, e_r\}) \to \\
H_0(G_2) \oplus H_0(G_1\setminus \{e_1,\ldots, e_r\})
\to H_0(G_1\setminus e_1) \to 0.
}
We have inductively $H_0(G_2\setminus \{e_2,\ldots,e_r\}) \cong H_0(G_2)$ and we
deduce
\eq{6.16}{H_0(G_1\setminus \{e_1,\ldots,e_r\}) \cong H_0(G_1\setminus e_1)  \cong H_0(G_1).
}
Let $f$ be as in \eqref{6.13} and assume there exists $c_{I_q,\dotsc,I_r}
\neq 0$ as above. We claim that $Y\cap E_1\cap\ldots\cap E_r$ can be
computed as follows. For clarity, it is convenient to change notation
a bit and write $D_i \subset P_i$ for the exceptional divisor. Abusively,
$E_i$ will denote any pullback of $D_i$ to a $P_j$ for $j>i$. Take the
strict transform $Y_1$ to $P_1$ and intersect with $D_1$. Now take the
strict transform ${\rm st}_{2,1}(Y_1\cap D_1)$ of $Y_1\cap D_1$ to $P_2$ and
intersect with $D_2$. continue in this fashion. The assertion is
\ml{6.17}{Y\cap \bigcap_1^r E_i = \\
E_r\cap st_{r,r-1}\Big(D_{r-1}\cap
st_{r-1,r-2}(D_{r-2}\cap\ldots st_{2,1}(D_1\cap Y_1)\ldots )\Big)
}
This is just an elaboration on \eqref{6.11}, \eqref{6.13}. The left hand side
amounts to taking the terms with $|I_1|+\ldots + |I_r|$ minimal, removing
appropriate powers of defining equations for the exceptional divisors, and
then restricting; while the right hand side takes those terms with
$\sum_{q}^r |I_j|$ minimum for $q=1,\dotsc,r-1$. By what we have seen,
these yield the same answer.

It remains to see that the intersection \eqref{6.17} doesn't contain the
vertex $v$. We have seen \eqref{6.7} that $D_1\cap Y_1$ is a union of the
pullbacks of graph hypersurfaces for $G_1$ and $\Gamma//G_1$. We have a
cartesian diagram \minCDarrowwidth.5cm
\eq{}{\begin{CD}E_1\cap D_2 \cong L_1\times \P^{p_2-1}\times
\P^{p_1-p_2-1} @>>>  P_2 @>>> BL(\lambda_2
\subset
\P^{p_1-1}) @>>>
\P^{p_2-1}\\ @VVV @VVV  @VVV\\
D_1\cong L_1\times \P^{p_1-1} @>>> P_1 @>\rho_1>> \P^{p_1-1}
\end{CD},
}
where $\lambda_2\cong \P^{p_2-1}$ corresponds to $L_2 \supset L_1$, and
the strict transform in $P_1$ is the pullback $\widetilde L_2=
\rho_1^{-1}(\lambda_2)$. Of course the picture continues in this fashion
all the way up. In the end, we get
\eq{}{L_1\times \P^{p_r-1}\times\P^{p_{r-1}-p_r} \times
\ldots \times \P^{p_1-p_2-1}.
}
The strict transform of $X$ here, by \eqref{6.17}, is the union of
pullbacks of graph hypersurfaces
\eq{18}{pr_{L_1}^{-1}X_{\Gamma//G_1} \cup pr_r^{-1}X_{G_r}\cup
pr_{r-1}^{-1}X_{G_{r-1}//G_r}\cup \ldots \cup pr_1^{-1}X_{G_1//G_2}.
}
Now each of the graphs involved has $h_1 = 1$, so each of the graph
hypersurfaces is linear. As we have seen, they involve all the edge
variables so they do not vanish at any of the vertices. This completes
the proof of Proposition \ref{prop6.3}(ii). Finally, the proof of (iii) is straightforward
from (ii). One uses the existence of local coordinates as in \eqref{6.13} with respect to
which the defining equation of the strict transform is a sum of monomials with
coefficients $>0$, and elements in the strict transform $\tilde\sigma$ of
$\sigma^{2n-1}(\R)$ have coordinates $\ge 0$. (Points in $Y\cap \tilde\sigma$ could be
specialized to vertices.)\end{proof}

We are now in a position to explicit the motive \eqref{0.1} associated to a primitive divergent graph
$\Gamma \subset \P^{2n-1}$. Let $P \xrightarrow{\pi} \P^{2n-1}$ be as in Proposition \ref{prop6.3}. Let
$\Delta \subset \P^{2n-1}$ be the union of the $2n$ coordinate hyperplanes. Let $B := \pi^*\Delta$ and
let $Y \subset P$ be the strict transform of the graph hypersurface $X=X_\Gamma$. Consider the motive
\eqref{0.1}
\eq{7.21}{H := H^{2n-1}(P\setminus Y,B\setminus B\cap Y).
}
By construction, 
\begin{prop}The divisor $B \subset P$ has normal crossings. The Hodge
  structure on the Betti realization $H_B$ has the following 
properties: \newline\noindent
(i) $H_B$ has weights in $[0,4n-2]$. $W_0H_B \cong \Q(0)$. \newline\noindent
(ii) The strict transform $\tilde\sigma$ of the chain $\sigma^{2n-1}(\R)$ in Proposition
\ref{prop6.3}(iii) represents an homology class in $H_{2n-1}(P\setminus Y,B\setminus B\cap Y)$. The
composition 
$$W_0H_B \inj H_B \xrightarrow{\int_{\tilde\sigma}} \Q
$$
is a vector space isomorphism. 
\end{prop}
\begin{proof} we have the exact sequence
\eq{}{0 \to H^{2n-2}(B\setminus Y\cap B)/ H^{2n-2}(P\setminus Y) \to H \to
  H^{2n-1}(P-Y) 
}
Write $B = \bigcup B_i,\ B^{(r)} = \coprod B_{i_1}\cap \ldots \cap
  B_{i_r}$. We have a spectral sequence of Hodge structures
\eq{}{E_1^{p,q} = H^q(B^{(p+1)}\setminus B^{(p+1)}\cap Y) \Rightarrow
  H^{p+q}(B\setminus B\cap Y)
}
From known properties of weights for open smooth varieties, we get an exact
  sequence
\eq{}{H^0(B^{(2n-2)}) \to H^0(B^{(2n-1)}) \to W_0H \to 0
}
An analogous calculation with $B$ replaced by $\Delta \subset \P^{2n-1}$
  yields $\Q(0)$ as cokernel. It is easy to see that blowing up strict
  transforms of linear spaces doesn't change this cokernel. This proves
  (i). Assertion (ii) is straightforward.
\end{proof}

An optimist might hope for a bit more. Whether for all primitive divergent graphs, or for an identifiable
subset of them, one would like that the maximal weight piece of $H_B$ should be Tate
\eq{7.22}{gr^W_{max}H_B = \Q(-p)^{\oplus r}. 
}
Further one would like that there should be a rank $1$ sub-Hodge structure $\iota: \Q(-p) \inj
gr^W_{max}H_B$ such that the image of $\eta_\Gamma \in H_{DR}$ in $gr^W_{max}H_{DR}$ spans
$\iota(\Q(-p))_{DR}$. Our main result is that this is true for wheel and spoke graphs, (sections
\ref{sect_ws}, \ref{sec12}).

\section{the motive II}\label{sect_motII}

In this section we consider the class of the graph hypersurface $[X_\Gamma]$ in the
Grothendieck group $K_{{\rm mot}}$ of quasi-projective varieties over $k$ with the relation $[X]
= [Y] + [X\setminus Y]$ for $Y$ closed in $X$. We assume $\Gamma$ has $N$ edges
and $n$ loops. The basic result of \cite{BB} is that $[X_\Gamma]$ can be quite
general. In particular, the motive of $X_\Gamma$ is not in general mixed Tate.

From the physicists' point
of view, of course, one is primarily interested in the period \eqref{5.10}. Results
in \cite{BB} do not exclude the possibility of some mixed Tate submotive yielding
this period. The methods of \cite{BB} seem to require graphs with physically unrealistic numbers of
edges, so it is worth looking more closely at
$[X_\Gamma]$. In this section we pursue a naive projection technique based on
the fact that graph and related polynomials have degree $\le 1$ in each variable.
We stratify
$X_\Gamma$ and examine whether the strata are mixed Tate. For $N = 2n \ge 12$,
we identify a possible non-mixed Tate stratum. 

Curiously,
the stratum we consider turns out to be mixed Tate in ``most'' cases, but with a
computer it is not difficult to generate cases where it may not be. We give such an
example with $12$ edges. Note however that Stembridge \cite{S} has shown that all graphs with $\le 12$
edges are mixed Tate, so the particular example we give must in fact be mixed Tate. Techniques 
and results in this section should be compared with \cite{S}, which predates our work.

The basic observation of Kontsevich is that for $X_\Gamma$
mixed Tate, there will exist a polynomial $P_\Gamma$ with $\Z$-coefficients such
that for any finite field $\F_q$ we have $\# X(\F_q) = P_\Gamma (q)$. Stembridge has implimented a
computer algorithm for checking this. It might be of interest to try some of our examples to see if they
satisfy Kontsevich's condition.

If we fix an edge $e$, by \eqref{3.10} we can write the graph polynomial
\eq{}{\Psi_\Gamma = A_e\cdot \Psi_{\Gamma\setminus e} + \Psi_{\Gamma/e}.
}
Projecting from the point $v_e$ defined by $A_e(v_e)=1,\ A_{e'}(v_e) = 0, e' \neq e$
yields $pr_e : \P^{N_1} \setminus \{v_e\} \to \P^{N-2}$ and
\eq{7.2}{X_\Gamma \setminus  pr_e^{-1}(X_{\Gamma\setminus e})\cap X_\Gamma \xrightarrow{\cong} \P^{N-2} \setminus
X_{\Gamma\setminus e}.
}
One might hope to stratify $X_\Gamma$ and try to analyse its motive in this way. We know,
however, by \cite{BB} that in general this motive is very rich, and such elementary
techniques will not suffice to understand it. Indeed, we have
\eq{}{pr_e^{-1}(X_{\Gamma\setminus e})\cap X_\Gamma = pr_e^{-1}(X_{\Gamma\setminus e}\cap X_{\Gamma/e}),
}
so already at the second step we must analyse an intersection of two graph hypersurfaces.
What is amusing is that, in fact, one can continue a bit further, and the process gives
some indication of where motivic complications might first arise.

\begin{lem}\label{lem7.1} Assume $\Gamma$ has $n$ loops and $2n$ edges. Enumerate the edge
variables
$A_1,\dotsc,A_{2n}$ in such a way that $A_1A_2\cdots A_n$ appears with coefficient $1$ in
$\Psi_\Gamma$. Then we can write
\eq{7.4}{\Psi_\Gamma = \det(m_{ij} + \delta_{ij}A_i)_{1\le i, j\le n};\quad m_{ij} =
m_{ij}(A_{n+1},\dotsc,A_{2n}).
}
In other words, the first $n$ variables appear only on the diagonal.
\end{lem}
\begin{proof}Let $T \subset \Gamma$ be the subgraph with edges $e_{n+1},\dotsc,e_{2n}$. Our
assumption implies that $T$ is a spanning tree, so $\Z[E_\Gamma] \cong H_1(\Gamma)\oplus
\Z e_{n+1} \oplus \ldots \oplus \Z e_{2n}$. the linear functionals $e_i^\vee$ thus induce
an isomorphism
\eq{}{(e_1^\vee,\dotsc,e_n^\vee): H_1(\Gamma) \cong \Z^n
}
With respect to this basis of $H_1(\Gamma)$ the rank $1$ quadratic forms $(e_i^\vee)^2$
correspond to the matrices with $1$ in position $(i,i)$ and zeroes elsewhere, for $1\le
i\le n$. Define $(m_{ij})$ to be the symmetric matrix associated to the quadratic form
$\sum_{n+1}^{2n} A_i(e_i^\vee)^2$. The assertion of the lemma is now clear.
\end{proof}

\begin{lem}[\cite{D}] \label{lem7.2} Let $\psi = \det(m_{ij} + \delta_{ij}A_i)_{1\le i,j \le n}$,
where the $m_{ij}$ are independent of $A_1,\dotsc,A_n$. For $1\le k \le n$ write $\psi^k :=
\frac{\partial}{\partial A_k}\psi$ and $\psi_k := \psi|_{A_k=0}$. For $I, J \subset
\{1,\dotsc,n\}$ with $\# I = \# J$, define $\psi(I,J)$ to be the determinant as above with
the rows in $I$ and the columns in $J$ removed. Let $1\le k, \ell \le n$ be distinct
integers and assume $k, \ell \not\in I\cup J$. Then
\ml{}{\psi(I,J)^{k\ell}\psi(I,J)_{kl} - \psi(I,J)_{k}^\ell\psi(I,J)_{l}^k = \\
\pm\psi(I\cup \{k\},J\cup \{\ell\}) \psi(I\cup \{\ell\},J\cup \{k\})
}
The two factors on the right have degrees $\le 1$ in $A_i$ for $i\le n$.
\end{lem}
\begin{proof}We can drop the rows in $I$ and the columns in $J$ to begin with and ignore
the $A_\nu$ for $\nu \not\in \{k,\ell\}$. In this way, we reduce to the following
assertion. Let $M$ be an $n\times n$ matrix with coefficients in a commutative ring.
Assume $n\ge 2$. Write $M(S,T)$ for the matrix with rows in $S$
and columns in $T$ deleted. Then
\ml{}{\det M(\{1,2\},\{1,2\})\cdot\det M - \det M(\{1\},\{1\})\cdot \det M(\{2\},\{2\})
= \\
-\det M(\{1\},\{2\})\cdot \det M(\{2\},\{1\})
}
(By convention, the determinant of a $0\times 0$-matrix is $1$.). This is a
straightforward exercise.
\end{proof}

We attempt to stratify our graph hypersurface $X_\Gamma$ using the above lemmas. To fix
ideas, we assume $\Gamma$ has $2n$ edges and $n$ loops.
\newline\newline\noindent Step 1.  We order the edges so $\Psi_\Gamma$ admits a
description as in Lemma \ref{lem7.1}. \newline\newline\noindent
Step 2. Project as in \eqref{7.2} with $e=e_1$, to conclude
\ml{}{[X_\Gamma] = [\P^{2n-2}] + [{\rm Cone}(X_{\Gamma\setminus e_1}\cap X_{\Gamma/e_1})] -
[X_{\Gamma\setminus e_1}\cap X_{\Gamma/e_1}] = \\
[\P^{2n-2}] + 1 + ([\A^1]-1)[X_{\Gamma\setminus e_1}\cap
X_{\Gamma/e_1}]
}
\newline\newline\noindent
Step 3. Using \eqref{3.10}, we can write (with notation as in Lemma \ref{lem7.2} and
$\Psi = \Psi_\Gamma$)
\ga{}{\Psi_{\Gamma \setminus e_1} = \frac{\partial}{\partial A_1}\Psi_\Gamma = A_2\Psi_{\Gamma \setminus \{
e_1, e_2\}} + \Psi_{(\Gamma \setminus  e_1)/e_2} = A_2\Psi^{12}+\Psi^1_2 \\
\Psi_{\Gamma/e_1} = \Psi_\Gamma|_{A_1=0} = A_2\Psi_{(\Gamma/e_1) \setminus e_2} +
\Psi_{\Gamma/\{e_1,e_2\}} = A_2\Psi^2_1 + \Psi_{12}. \notag
}
Eliminating $A_2$, we conclude that projection from $\P^{2n-2}$ onto $\P^{2n-3}$ with
coordinates
$A_3,\dotsc,A_{2n}$ carries $X_{\Gamma-e_1}\cap X_{\Gamma/e_1}$ onto the hypersurface
defined by $\Psi^1_2\Psi^2_1 - \Psi^{12}\Psi_{12}=0$. By Lemma \ref{lem7.2},
\eq{7.10}{\Psi^1_2\Psi^2_1 - \Psi^{12}\Psi_{12}= \Psi(1,2)\Psi(2,1) = \Psi(1,2)^2
}
(The right hand identity holds because $\Psi = \Psi_\Gamma$ is the determinant of a
symmetric matrix.) \newline\newline\noindent
Step 4. Write $\sV(I)$ for the locus of zeroes of a homogeneous ideal $I$. The projection
in step 3 blows up on $\sV(\Psi^1_2,\Psi^2_1,\Psi^{12},\Psi_{12})$, and we conclude
\ml{7.11}{[X_{\Gamma\setminus e_1}\cap X_{\Gamma/e_1}] = \\ [X(1,2)]+[{\rm Cone}
\ \sV(\Psi^1_2,\Psi^2_1,\Psi^{12},\Psi_{12})] -
[\sV(\Psi^1_2,\Psi^2_1,\Psi^{12},\Psi_{12})]  = \\
[X(1,2)]+1+([\A^1]-1)[\sV(\Psi^1_2,\Psi^2_1,\Psi^{12},\Psi_{12})].
}
\newline\newline\noindent
Step 5. One could try to study the motive of $\sV(\Psi^1_2,\Psi^2_1,\Psi^{12},\Psi_{12})$,
but the elimination theory gets complicated, so instead we focus on $[X(1,2)]$. Since
$\Psi(1,2)$ has degree $\le 1$ in $A_3$ we may project onto $\P^{2n-4}$ with coordinates
$A_4,\dotsc,A_{2n}$. It might seem that we could repeat the argument starting from step 2
above, but there is a problem. Writing $\Psi = \det M$ with $M$ symmetric, we have
$\Psi(1,2) = \det M(1,2)$, where $M(1,2)$ is obtained from $M$ by deleting the first row
and the second column. This matrix is no longer symmetric. Just as in \eqref{7.2}, the
projection $X(1,2) \to \P^{2n-4}$ blows up over $\sV(\Psi(1,2)^3, \Psi(1,2)_3)$.
\newline\newline\noindent
Step 6. Just as in step 3, we project $\sV(\Psi(1,2)^3, \Psi(1,2)_3)$ to $\P^{2n-5}$ with
coordinates $A_5,\dotsc,A_{2n}$. When we eliminate $A_4$ we find the image of the
projection is given by the zeroes of
\ml{7.12}{\Psi(1,2)^{34}\Psi(1,2)_{34} - \Psi(1,2)^{3}_4\Psi(1,2)^{4}_3 \stackrel{{\rm Lemma} \
\ref{lem7.2}}{=} \\
\Psi(\{1,3\},\{2,4\})\cdot \Psi(\{1,4\},\{2,3\}).
}
\newline\newline\noindent
Step 7. At this point something new has happened. The right hand side in \eqref{7.12} is
not a square. Although both factors have degree $\le 1$ in $A_5$, we will at the next
stage in our motivic stratification have to deal with
\eq{7.13}{\sV(\Psi(\{1,3\},\{2,4\}), \Psi(\{1,4\},\{2,3\})).
}
Here Lemma \ref{lem7.2} no longer applies. We find by example that eliminating $A_5$, the
resulting hypersurface in $\P^{2n-6}$ in general no longer factors into factors with
degrees $\le 1$ in $A_6$. Projection then is no longer an isomorphism at the generic
point, and the argument is blocked.

\begin{ex}The computer yields the following example of a graph with 6 loops and 12 edges
for which the projection \eqref{7.13} has an irreducible factor with degree $2$ in $A_6$.
Take $7$ vertices labeled $1,2,\dots,7$ and connect them with edges as indicated:
\ml{}{(1,2), (2,3), (3,4), (4,5), (5,6), (6,7), (7,2), (7,3), \\
(6,4), (5,1), (5,3), (4,1)
}
\end{ex}
Note that this graph is mixed Tate though by explicit computation,
which finds it $\sim \zeta(3)\zeta(5)$.
\section{General Remarks}\label{sec9}

Let $\Gamma$ be a graph with $n$ loops and $2n$ edges. We assume all subgraphs of $\Gamma$
are convergent so the period $P(\Gamma)$ is defined (Proposition \ref{prop4.2}). The
Schwinger trick (Corollary \ref{cor5.3}) relates $P(\Gamma)$ to an integral computed in
Schwinger coordinates in $\P^{2n-1}$. To avoid confusion, we write $P_{{\rm quadric}}(\Gamma)$ for
the period \eqref{4.7} of the configuration of Feynman quadrics associated to
$\Gamma$ and $P_{{\rm graph}}(\Gamma)$ for the graph period. We have by \eqref{5.10}
\eq{}{P_{{\rm quadric}}(\Gamma) \in \Q^\times \pi^{-2n}P_{{\rm graph}}(\Gamma).
}
Proposition
\ref{prop6.3} shows that there is a suitable birational transformation   $\pi: P\to
\P^{2n-1}$ defined over $\Q$, such that the integrand
$\eta\in \Gamma(\P^{2n-1}, \omega(2X))$
keeps poles only along the strict transform $Y$ of the
discriminant hypersurface $X$, that is $\pi^*(\eta)\in \Gamma(P, \omega(2Y))$. Thus,
denoting by $B$ the total transform of the union $\Delta$ of coordinate hyperplanes $A_i=0$,
the form $\eta$ yields a class
\ga{9.1}{\pi^*\eta\in \Gamma(P, \omega(2Y))\to H^{2n-1}_{DR}(P\setminus Y, B\setminus B\cap Y)}
in relative de Rham cohomology. On the other hand, Proposition \ref{prop6.3} shows
that the strict transform $\tilde{\sigma}^{ 2n-1}(\R)$ of the cycle of integation
is disjoint from $Y$. Thus it yields a relative homology class
\ga{9.2}{\tilde{\sigma}^{ 2n-1}(\R)  \in H_{2n-1}(P\setminus Y, B\setminus B\cap Y)=
H^{2n-1}_{{\rm Betti}}(P\setminus Y, B\setminus B\cap Y)^\vee}
in Betti cohomology. More precisely
\begin{claim} \label{claim9.1} The period integral \eqref{4.3} $P_{{\rm quadric}}(\Gamma) \in
\pi^{-2n}\Q^\times\cdot P_{{\rm graph}}(\Gamma)$, where $P_{{\rm graph}}(\Gamma)$ is a period of the
cohomology
$H^{2n-1}(P\setminus Y, B\setminus B\cap Y)$. By period here we mean the integral of an
algebraic de Rham form $\pi^*\eta$ defined over $\Q$ against a $\Q$-homology chain
$\tilde{\sigma}^{ 2n-1}$.
\end{claim}
Suppose now, as has been established in a number of cases \cite{BK}, that the period is
related to a zeta value:
$P_{{\rm quadric}}(\Gamma) \in \pi^{\Z}\Q^\times \zeta(p)$. Then the general guideline for what we
wish to understand is the following.

One has now a good candidate for a triangulated category of mixed motives
over $\Q$, defined by Voevodsky, Levine and Hanamura 
(\cite{DG}, section 1  and references there for the discussion here). One
further considers the triangulated subcategory spanned by 
$\Q(n)$, $n\in \Z$. In this category, one has
\ga{9.3}{{\rm Hom}^j(\Q(0), \Q(p))=\begin{cases}
\Q & p=j=0\\
K_{2p-1}(\Q)\otimes \Q & p\ge 1, \ j=1\\
0 & {\rm else}
  \end{cases}}
The iterated extensions of $\Q(n)$ form an abelian subcategory which is the heart of a $t$-structure.

Borel's work on the $K$-theory of number fields \cite{Bo}, \cite{So} tells us that 
$K_{2p-1}(\Q)\otimes \Q \cong \Q$ for $p=2n-3, n\ge 2$, so there is a one dimensional space of motivic extensions of
$\Q(0)$ by $\Q(p)$. We want to understand their periods. Let $E$ be a
nontrivial such extension.
We write $E_{DR} = \Q\cdot e_0 \oplus \Q\cdot
e_p$, with $F^0E_{DR} = \Q e_0$. The Betti realization is $E_\C = \C\cdot e_0 \oplus \C\cdot
e_p$ and $E_\Q = \Q\cdot (2\pi i)^pe_p \oplus \Q\cdot(e_0 + \beta e_p)$ for a
suitable $\beta$. The corresponding Hodge structures on the $\Q(i)$ are
\eq{}{(\Q(0)_{DR} = \Q\cdot\epsilon_0, \Q(0)_\Q = \Q\cdot \epsilon_0),\quad
  (\Q(p)_{DR} = \Q\cdot\epsilon_p, \Q(p)_\Q = \Q\cdot (2\pi i)^p\epsilon_p)
}
We have an exact sequence
\eq{}{0 \to \Q(p) \to E \to \Q(0) \to 0
}
given by $\epsilon_p \mapsto e_p,\ e_0 \mapsto \epsilon_0$. The ambiguity
here is that we can replace $e_0+ \beta e_p$ by $e_0 + (\beta + c(2\pi i)^p) e_p$ for $c\in
\Q$ as a basis element for $E_\Q$, so $\beta \in \C/(2\pi i)^p\Q$ is well
defined. In fact, $\text{Ext}^1_{MHS}(\Q(0),\Q(p)) = \C/(2\pi i)^p\Q$ and
$\beta$ is the class of $E$.

To compute the period, consider the dual object $E^\vee$, with $E^\vee_{DR}
= \Q e_0^\vee \oplus \Q e_p^\vee$ and $E^\vee_\Q = \Q e_0^\vee \oplus
\Q(2\pi i)^{-p}(e_p^\vee - \beta e_0^\vee)$. By definition, the period is
obtained by pairing $F^0E_{DR}$ against a lifting in $E^\vee_\Q$ of the
generator $(2\pi i)^pe_p^\vee \in \Q(-p)_\Q = E^\vee_\Q/\Q(0)_\Q$. This
yields
\eq{}{ \langle e_0, (2\pi i)^{-p}(e_p^\vee - \beta e_0^\vee) \rangle = -(2\pi
  i)^{-p}\beta.
}
It is better from the period viewpoint to dualize and consider the
period of $E^\vee$, which is an extension of $\Q(-p)$ by $\Q(0)$. this
yields
\eq{}{\langle e_p^\vee,e_0+ \beta e_p\rangle = \beta
}

For $E$ a non-split motivic extension of $\Q(0)$ by $\Q(p)$, $p$ odd, $\ge
3$, let $\beta \in \C/(2\pi i)^p\Q$ be the extension class. Note
$\text{Im}(\beta) \in \R$ is well defined. one knows by
the Borel regulator theory \cite{Bo}, \cite{So} that $\zeta(p) \in \text{Im}(\beta)\Q^\times$.

Now consider our graph $\Gamma$ with period related to $\zeta(p)$. The motive
$H^{2n-1}(P\setminus Y, B\setminus B\cap Y)$ has lowest weight piece $\Q(0)$, so we might expect to find
inside it a subquotient motive of rank $2$ which is an extension of $\Q(-p)$ by $\Q(0)$. By the above discussion,
we would then hope
\eq{}{P_{{\rm graph}}(\Gamma) \in \zeta(p)\Q^\times.
}
By \eqref{5.10} this would yield $P_{{\rm quadric}}(\Gamma) \in \pi^{-2n}\zeta(p)\Q^\times$. For example,
take $\Gamma = \Gamma_n$ to be the wheel with $n$ spokes. Then $p=2n-3$ and we expect, if indeed the
$\zeta$-values computed in \cite{BK} are motivic, to find
\eq{}{P_{{\rm graph}}(\Gamma_n) \in \zeta(2n-3)\Q^\times;\quad P_{{\rm quadric}}(\Gamma_n) \in
\pi^{-2n}\zeta(2n-3)\Q^\times }

 The aim of the next sections is to show for the wheel and spoke family of
examples what can be done motivically. We will show in particular
\ga{9.5}{H^{2n-1}(\P^{2n-1}\setminus X)=\Q(-2n+3).}
Moreover, $H^{2n-1}_{DR}(\P^{2n-1}\setminus X)$ is spanned by $\eta$. Even in this special
case, we are not able to find a suitable rank $2$ subquotient motive of $H^{2n-1}(P\setminus Y,
B\setminus B\cap Y)$.

\section{Correspondences} \label{sec10}
We will assume in this section that $\Gamma$ has $n$ loops and $2n$ edges.
So one has $2n$ Feynman quadrics which we denote by $q_e$, of equation $Q_e$, see section 5. Recall concretely that to an edge $e$, one associates coordinates $x_e(i), i=1,\ldots, 4=j$. Given  an orientation of $\Gamma$, to  a vertex $v$, one associates the relation $\sum_e {\rm sign}(v,e) x_e(i)=0$ for all $i=1,\ldots, j=4$. Then $q_e=:q_e^j$ is defined by $Q_e^j:=\sum_{a=1}^{4=j}
x_e(a)^2=0$ in $\P^{jn-1}$. One defines $\sQ=\sQ_j\subset \P^{jn-1}\times \P^{2n-1}$ by the equation
$\sum_e A_eQ_e^j=0$. This defines a correspondence
\ga{10.1}{\begin{CD}
\P^{2n-1}\times \P^{jn-1}\setminus \sQ_j @>\A^{2n-1}-{\rm fibration}>> \P^{jn-1}\setminus \cap_{e=1}^{2n} q_e^j\\
@V\pi_jVV\\
\P^{2n-1}
\end{CD}}
We discuss now this correspondence for the Feynman quadrics, i.e. $j=4$.
 On the other hand, we can consider all the definitions above for other $j$, and we discuss the resulting correspondence \eqref{10.1} for $j=1$ and $j=2$ as well.

For $j=1$, we rather consider the projection ${\rm proj}: \sQ_1\to \P^{2n-1}$. Let us denote by $\Sigma\subset \Q_1$ the closed subscheme with ${\rm proj}^{-1}(x)\cap \Sigma={\rm Sing}({\rm proj}^{-1}(x)).$ Then $\Sigma \to X$ is the desingularization $\P(N)\to X$ studied in Proposition
\ref{prop44.2}.

We assume now $j=2$. Recall that if $Z\subset \P^{2N+1}$ is a smooth even dimensional quadric, then
\ga{10.2}{
H^{j}_c(\P^{2N+1}\setminus Z)=
\begin{cases}
0 & j\neq 2N\\
\Q(-N)[\ell_1-\ell_2] & j=2N
 \end{cases}}
where $\ell_i$ are the 2 rulings of $ Z$.
We define
\ga{10.3}{X_i=\{(A)\in \P^{2n-1}, {\rm rk} (\sum_e A_e Q_e^1)< n-i\}.}
So $X=X_0$, and $X_{i+1}$ is the singular locus of $X_i$. We denote by $j=j_0: \P^{2n-1}\setminus X \to \P^{2n-1}, j_i: X_{i-1}\setminus X_i\to X_{i-1}$.
Over $X_i$, the quadric $\sum_eA_eq_e^j$ is a cone over a smooth quadric
$\overline{ \sum_eA_eq_e^j} \subset \P^{j(n-i)-1}$, thus by homotopy invariance
and base change for $R(\pi_j)_!$ (\cite{DeSGA}), one  obtains
\begin{prop} \label{prop10.1}
\ga{10.4}{
R^i(\pi_4)_!\Q=\begin{cases}
j_! \Q(-2n+1) & i=4n-1\\
(j_1)_! \Q(-2n-1) & i=4n+3\\
\ldots & \ldots\\
(j_a)_! \Q(-2n+1-2a) & i=4n+4a
\end{cases}}
\ga{10.5}{
R^i(\pi_2)_!\Q=\begin{cases}
j_! \Q(-n+1) & i=2n-1\\
(j_1)_! \Q(-n) & i=4n+1\\
\ldots & \ldots\\
(j_a)_! \Q(-2n+1-a) & i=2n+2a
\end{cases}
}
\end{prop}
We draw now two consequences from this computation.
\begin{prop} \label{prop10.2}
One has maps
\ga{10.6}{
H^{2n-1}_c(\P^{2n-1}\setminus X)\to H^{2n}_c(\P^{4n-1}\setminus \cap_{e=1}^{2n} q_e^4) \to
H^{4n-1}_c(\P^{4n-1}\setminus \cup_{e=1}^{2n} q_e^4)\\
{\rm in \ particular \ dually} \notag \\
H^{4n-1}(\P^{4n-1}\setminus \cup_{e=1}^{2n}q_e^4)(2n)\to H^{2n-1}(\P^{2n-1}\setminus X).\notag}
\end{prop}
\begin{proof}
By \eqref{10.4}, the term $E_2^{2n-1, 4n-1}=H^{2n-1}_c(\P^{2n-1}\setminus X)(-2n+1)$ of the Leray spectral sequence for $\pi_4$ maps to $H^{2n-1+4n-1}_c(\P^{2n-1}\times \P^{4n-1}\setminus \sQ_4)$, which in turns is equal to $H^{2n}_c(\P^{4n-1}\setminus \cap_{e=1}^{2n} q_e^4)(-2n+1) $ by homotopy invariance. The second map comes from the Mayer-Vietoris spectral sequence for
$\cup_{e=1}^{2n} q_e^4$. \end{proof}
\begin{rmk} \label{rmk10.3}
We will see in  section 11 on the wheel with $n$ spokes that for $n=3$, the first map is an isomorphism, but in genaral, we do not control it. \end{rmk}

\begin{prop} \label{prop10.4}
Assume $\cap_{e=1}^{2n}q_e^2\neq \emptyset$, for example for the wheel with $n$ spokes (see section 11). Then $$H^{2n-1}_c(\P^{2n-1}\setminus X)=H^{2n-2}(X)/H^{2n-2}(\P^{2n-1})$$ is supported along $X_a$ for some $a\ge 1$. \end{prop}
\begin{proof}
By homotopy invariance again and by assumption, we have
\ga{10.7}{
H^{2n-1+2n-1}_c(\P^{2n-1}\times \P^{2n-1}\setminus \sQ_2)=H^0_c(\P^{2n-1}\setminus \cap_{e=1}^{2n}q_e^2)=0.}
So the Leray spectral sequence for $\pi_2$ together with
\eqref{10.4} imply that $E_\infty^{2n-1, 2n-1}=0$, with $E_2^{2n-1,2n-1}=H^{2n-1}_c(\P^{2n-1}\setminus X)(-n+1)$. So, since $R^i(\pi_2)_!$ is supported in lower strata of $X$, this shows the proposition.

\end{proof}
\begin{rmk} \label{rmk10.5}
We will see in  section 11 on the wheel with $n$ spokes that for $n=3$, the Leray spectral sequence will equate $H^0(X_1)(-1)\xrightarrow{\cong} H^4(X)/H^4(\P^5)$.
\end{rmk}

\section{Wheel and Spokes}\label{sect_ws}

The purpose of this section is to compute the middle dimensional
cohomology for a graph polynomial in a non-trivial case.
The geometry we will be using involves only projections,  homotopy invariance and Artin vanishing theorem. Consequently, our cohomology computation holds for Betti or \'etale cohomology, and would for motivic cohomology if one had Artin vanishing. To unify notations, we
denote this cohomology as $H(?, \Q)$ rather than $\Q_\ell$ in the $\ell$-adic case.

Fix $n \ge 3$
and let $\Gamma = WS_n$ be the graph which is a wheel with $n$ spokes.
$WS_n$ has vertices $\{0, 1, \dotsc,n\}$ and edges $e_i = (0,i),\ 1\le
i\le n$ and $e_j = (j-n,j-n +1 \mod n),\ n+1 \le j \le 2n$. Suitably
oriented, $\ell_i = e_i + e_{i+n} - e_{i+1\mod n},\ 1\le i\le n$ form a
basis for the loops. The following is straightforward.

\begin{lem}$\Gamma$ has $n$ loops and $2n$ edges. Every proper subgraph
$\Gamma' \subsetneq \Gamma$ is convergent so the period
$P(\Gamma)$ is defined (see Proposition \ref{prop4.2}).
\label{nosubdiv}\end{lem}
\begin{proof}Omitted.  \end{proof}

Let $T_i, 1 \le i\le 2n$ be variables. The graph polynomial of $\Gamma$
can be written
\eq{}{\Psi_\Gamma(T) = \det(\sum_{i=1}^{2n} T_i M^{(i)})
}
where
\eq{}{M^{(i)} = (M^{(i)}_{pq})_{1\le p, q\le n};\quad M^{(i)}_{pq} =
e_i^\vee(\ell_p)e_i^\vee(\ell_q).
}
It follows easily that
\ml{}{\Psi_\Gamma = \\
\det \begin{pmatrix} T_1+T_2+T_{n+1} &
-T_2 & 0 & \hdots & 0 & -T_1 \\
-T_2 & T_2+T_3+T_{n+2} & -T_3 & \hdots & 0 & 0 \\
\vdots & \vdots & \vdots & \hdots & \vdots & \vdots \\
-T_1 & 0 & 0 & \hdots & -T_n & T_n + T_1 + T_{2n}\end{pmatrix}.
}
It will be convenient to make the change of variables
\ga{11.4}{B_i = T_{i+1} + T_{i+2} + T_{i+1+n}, \
A_i = -T_{i-2}
}
where all the indices are counted modulo $n$ and taken in $[0,\ldots, n]$.
Write
\eq{11.5}{\Psi_n = \Psi_n(A,B) = \det \begin{pmatrix}B_0 & A_0 & 0 &
\hdots & \hdots & A_{n-1} \\
A_0 & B_1 & A_1 & \hdots &\hdots & 0 \\
0 & A_1 & B_2 & A_3 &\hdots & 0 \\
\vdots & \vdots & \vdots & \vdots & \hdots & \vdots \\
A_{n-1} & 0 & \hdots & \hdots & A_{n-2} & B_{n-1}\end{pmatrix}.
}
The graph hypersurface in the $A, B$-coordinates is given by
\eq{}{\P^{2n-1}\supset X_n : \Psi_n(A,B) = 0.
}

Define $H^*(X_n,\Q)_{{\rm prim}} := \text{coker}(H^*(\P^{2n-1},\Q) \to
H^*(X_n,\Q))$. We formulate now our main theorem.
\begin{thm}\label{thm11.2}
Let $X_n\subset \P^{2n-1}$ be the graph polynomial hypersurface
for the wheel with $n\ge 3$ spokes. Then one has
$$H^{2n-1}(\P^{2n-1}\setminus X_n)\cong \Q(-2n+3)$$
or equivalently, via duality
$$H^{2n-2}(X_n,
\Q)_{{\rm prim}}
\cong \Q(-2).$$ In particular, $H^{2n-1}(X_n, \Q)_{{\rm prim}}$ is independent of $n \ge 3$.
\end{thm}
\begin{proof} The proof is quite long and involves several geometric steps.
We first define homogeneous polynomials $Q_{n-1}$ and $K_n$ as
indicated:
\ml{11.7}{\Psi_n=B_0Q_{n-1}(B_1,\ldots, B_{n-1}, A_1,\ldots, A_{n-2})+\\
K_n(B_1,\ldots, B_{n-1}, A_0,\ldots, A_{n-1})
}

Here
\ml{11.8}{Q_{n-1}(B_1,\ldots, B_{n-1}, A_1,\ldots, A_{n-2})= \\
{\rm det} \begin{pmatrix}
B_1 & A_1 & 0& \ldots & \ldots & 0\\
A_1 & B_2 & A_2 & 0 &\ldots & 0\\
\ldots & \ldots & \ldots & \ldots & \ldots \\
0&\ldots & \ldots & \ldots &A_{n-2} & B_{n-1}
\end{pmatrix}
}
\begin{lem}One has inductive formulae:
\ml{11.9}{
Q_{n-1}(B_1,\ldots, B_{n-1}, A_1,\ldots, A_{n-2})=\\
B_1Q_{n-2}(B_2,\ldots, B_{n-1}, A_2,\ldots, A_{n-2}) -\\
A_1^2
Q_{n-3}(B_3,\ldots, B_{n-1}, A_3,\ldots, A_{n-2})= \\
B_{n-1}Q_{n-2}(B_1,\ldots, B_{n-2}, A_1,\ldots, A_{n-3}) -\\
A_{n-2}^2
Q_{n-3}(B_1,\ldots, B_{n-3}, A_1,\ldots, A_{n-4});
}
and
\ml{11.10}{K_n(B_1,\ldots, B_{n-1}, A_0,\ldots, A_{n-1})=\\
-A_0^2Q_{n-2}(B_2,\ldots, B_{n-1}, A_2,\ldots, A_{n-2}) -\\
A_{n-1}^2 Q_{n-2}(B_1,\ldots, B_{n-2}, A_1,\ldots, A_{n-3})
+2(-1)^{n-1} A_0\cdots A_{n-1}.
}
\end{lem}
\begin{proof} Straightforward. \end{proof}
The following lemma is a direct application of Artin's vanishing theorem \cite{Ar}, Th\'eor\`eme 3.1,  and homotopy invariance, and will be the key ingredient to the computation.
\begin{lem} \label{lem11.4}
Let $V\subset \P^N$ be a hypersurface which is a cone over the
hypersurface  $W\subset \P^a$. Then one has
$$ H^i(\P^N\setminus V)=0  \ {\rm for} \
i>a$$ or equivalently
$$H^j_c(\P^N\setminus V)=0 \  {\rm for}  \ j<2N-a.$$
\end{lem}
\begin{proof} The projection  $\P^N\setminus V\to \P^a\setminus W$ is a $\A^{N-a}$-fibration
By homotopy invariance, $H^j_c(\P^N\setminus V)=H^{j-2(N-a)}_c(\P^a\setminus W)(-(N-a))$ and
by Artin's vanishing
$H^{j-2(N-a)}_c(\P^a\setminus W)=0$ for $j -2(N-a)<a$ i.e. for $j<2N-a$.
\end{proof}

For a homogeneous ideal $I$ or a finite set $F_1, F_2, \ldots$ of
homogeneous polynomials, we write $\sV(I)$ or $\sV(F_1, F_2, \ldots)$
for the corresponding projective scheme. We will need to pass back and
forth via various projections. In confusing situations we will try to
specify the ambiant projective space. A
superscript $(i)$ will mean the ambient projective space is $\P^i$. In
the following lemma,
$\P^{2n-1}$ has coordinates
$(B_0:\dotsc: B_{n-1}:A_0:\dotsc:A_{n-1})$ and $\P^{2n-2}$ drops the $B_0$.

\begin{lem}\label{lem11.5} We have
\eq{11.11}{H^{2n-2}(X_n,\Q) \cong
H^{2n-4}\Big(\sV(Q_{n-1},K_n)^{(2n-2)},\Q(-1)\Big).
}
\end{lem}
\begin{proof} By \eqref{11.7}, one has
\eq{11.12}{X_n \cap \sV(Q_{n-1}) = \sV(Q_{n-1},K_n)^{(2n-1)}.
}
Let $p=(1,0,\dotsc,0)\in
\P^{2n-1}$. Projection from $p$
gives an isomorphism (use \eqref{11.7} to solve for $B_0$)
\eq{11.13}{\pi_p: X_n \setminus  X_n\cap \sV(Q_{n-1}) \cong \P^{2n-2} \setminus
\sV(Q_{n-1}).  }
We get a long exact sequence
\ml{11.14}{H^{2n-2}_c(\P^{2n-2}\setminus\sV(Q_{n-1})) \to H^{2n-2}(X_n) \to \\
H^{2n-2}(\sV(K_n,Q_{n-1})^{(2n-1)})
\to  H^{2n-1}_c(\P^{2n-2}\setminus\sV(Q_{n-1}))
}
Since the polynomial $Q_{n-1}$ does not involve $A_0$ or $A_{n-1}$, we
can apply Lemma \ref{lem11.4} with $N=2n-2$ and $a=2n-4$ to deduce
\eq{11.15}{H^{i}_c(\P^{2n-2}\setminus \sV(Q_{n-1})) = (0),\quad i<2n.
}
We conclude
\eq{}{H^{2n-2}(X_n) \cong H^{2n-2}(\sV(K_n,Q_{n-1})^{(2n-1)}).
}
The projection $\pi_p$ is an $\A^1$-fibration
$$\sV(K_n, Q_{n-1})^{(2n-1)}-p
 \to \sV(K_n, Q_{n-1})^{(2n-2)}$$
and we obtain
\ml{}{H^{2n-2}(\sV(K_n, Q_{n-1})^{(2n-1)})\cong \ (2n-2>0)\\
  H^{2n-2}_c(\sV(K_n, Q_{n-1})^{(2n-1)}-p) \cong \\
H^{2n-4}(\sV(K_n, Q_{n-1})^{(2n-2)})(-1).
}
\end{proof}

We now consider the line $\ell$ with coordinate functions $A_0, A_{n-1}$
\ga{11.18}{\ell \subset \P^{2n-2}(B_1:\ldots:B_{n-1}:A_0:\ldots:A_{n-1})
\\
  \ell:B_1=\ldots = B_{n-1}=A_1=\ldots=A_{n-2} = 0.\notag
}
One has $\ell \subset \sV(Q_{n-1},K_n)^{(2n-2)}$. The sequence
\ml{11.19}{0 \to H^{2n-4}_c(\sV(Q_{n-1}, K_n)^{(2n-2)}\setminus \ell) \\
\to H^{2n-4}(\sV(Q_{n-1}, K_n)^{(2n-2)}) \to H^{2n-4}(\ell)
}
 together with the previous lemma implies
\eq{11.20}{H^{2n-4}_c(\sV(Q_{n-1}, K_n)^{(2n-2)}\setminus \ell)(-1)
\cong
H^{2n-2}( \tilde{X}_n,\Q)
}
where $$H^{2n-2}(X_n)=H^{2n-2}(\tilde{X}_n) \ {\rm for} \  n>3$$ and
for $n=3$,
$$H^4(\tilde{X}_3)={\rm ker}(H^4(X_3)\to H^2(\ell)(-1))\cong H^4(X_3)_{{\rm prim}}.$$
The next step is now motivated by the shape of the matrix \eqref{11.5}. If we wish to induct on $n$, we have to find the geometry which gets rid of the corner term $A_{n-1}$ in the matrix.
We project further to
$\P^{2n-4} = \P^{2n-4}(B_1:\ldots:B_{n-1}:A_1:\ldots:A_{n-2})$. Let
\eq{11.21}{r:\sV(Q_{n-1}, K_n)^{(2n-2)}\setminus \ell \to \sV(Q_{n-1})^{(2n-4)}
}
be the projection with center $\ell$. It is clear from \eqref{11.10} that
the fibres of $r$ are conics in the variables $A_0, A_{n-1}$ with
discriminant
\ml{11.22}{\delta_{n-1}(B_1,\ldots, B_{n-1}, A_1,\ldots, A_{n-2}) :=\\ Q_{n-2}(B_2,\dotsc,B_{n-1},A_2,\dotsc,A_{n-2})\cdot
  Q_{n-2}(B_1,\dotsc,B_{n-2},A_1,\dotsc,A_{n-3}) \\
  - (A_1\cdots A_{n-2})^2
}
We show that in fact the situation is degenerated:
\begin{lem}\label{lem11.6} One has
\ga{11.23}{\delta_{n-1}(B_1,\ldots, B_{n-1}, A_1,\ldots, A_{n-2})= \\
Q_{n-3}(B_2,\ldots, B_{n-2}, A_2,\ldots, B_{n-3})\cdot
Q_{n-1}(B_1,\ldots, B_{n-1}, A_1,\ldots, A_{n-2}).\notag
}
In particular, the general fibre of $r$ in \eqref{11.21} is a double line
(so
$\{Q_{n-1}=K_n=0\}$ is non-reduced)
\end{lem}
\begin{proof} We compute in the ring
\ga{11.24}{K[B_1,\ldots, B_{n-1}, A_1, \ldots, A_{n-2}, \frac{1}{
Q_{n-2}(B_2,\ldots, B_{n-1}, A_2,\ldots, A_{n-2}}].}
One has
\ml{11.25}{B_1=\\
A_1^2 Q_{n-3}(B_3,\ldots, B_{n-1},A_3,\dotsc,A_{n-2})/
Q_{n-2}(B_2,\ldots, B_{n-1},
A_2,\ldots, A_{n-2}) \\
+Q_{n-1}(B_1,\ldots, B_{n-1}, A_1,\ldots, A_{n-2})/Q_{n-2}(B_2,\ldots, B_{n-1},
A_2,\ldots, A_{n-2}).
 }
This yields
\ga{11.26}{\delta_{n-1}=\\
A_1^2\Big(\delta_{n-2}(B_2,\ldots, B_{n-1}, A_2, \ldots, A_{n-2})-\notag \\
Q_{n-2}(B_2,\ldots, B_{n-1}, A_2, \ldots, A_{n-2})\cdot
Q_{n-4}(B_3,\ldots, B_{n-2}, A_3,\ldots, A_{n-3})\Big)
 + \notag \\
Q_{n-1}(B_1,\ldots, B_{n-1}, A_1,\ldots, A_{n-2})\cdot Q_{n-3}(B_2,\ldots, B_{n-2}, A_2,\ldots A_{n-3}).\notag}
We now argue by induction starting with $n= 3$:
\ga{11.27}{\delta_{3-1}=B_1B_2-A_1^2=Q_2(B_1,B_2,A_1)\cdot 1.
}
\end{proof}
>From Lemma \ref{lem11.6} we see that the reduced scheme
$\sV(Q_{n-1}, K_n)_{{\rm red}} \setminus \ell$ is fibred over $\sV(Q_{n-1})^{(2n-4)}
\subset
\P^{2n-4}$ with general fibre $\A^1$. The fibres jump to $\A^2$ over the
closed set
\ml{11.28}{Z_{n-1}: \sV\Big(Q_{n-1}(B_1,\dotsc,B_{n-1},
A_1,\dotsc,A_{n-2}),
\\ Q_{n-2}(B_1,\dotsc,B_{n-2}A_1,\dotsc,A_{n-3}),
Q_{n-2}(B_2,\dotsc,B_{n-1},A_2,\dotsc, A_{n-2})\Big)
}
As a consequence, we
get an exact sequence
\ml{11.29}{H^{2n-9}(Z_{n-1})(-3) \to
  H^{2n-6}_c(\sV(Q_{n-1})^{(2n-4)}\setminus Z_{n-1})(-2) \to \\
  H^{2n-2}(\tilde {X}_n) \to H^{2n-8}(Z_{n-1})(-3)
}
with the tilde as in \eqref{11.20}.
\begin{lem}\label{lem11.7}
\begin{itemize}
\item[(i)] The restriction map $H^i(\P^{2n-4}) \to
  H^i(Z_{n-1})$ is
  surjective for $i < 2n-7$.
  \item[(ii)] $Z_2=\emptyset$.
  \item[(iii)] For $n\ge 4$ we have
  $$H^{2n-7}(Z_{n-1}) \cong H^{2n-6}_c(\{Q_{n-1}=0\}^{(2n-4)}\setminus Z_{n-1}).$$ \end{itemize}
\end{lem}
\begin{proof} (i) $Z_{n-1}$ is defined by 3 equations, thus by Artin's vanishing theorem
    $H^i_c(\P^{2n-4} \setminus Z_{n-1})=0$ vanishes for $i<2n-6$.\\
    (ii)  One has $Z_2 : B_1B_2-A_1^2=B_1=B_2=0$ in
  $\P^2(B_1:B_2:A_1)$, so $Z_2 = \emptyset$.\\
  (iii) For $n\ge 4$ we have
\ml{}{H^{2n-7}(\sV(Q_{n-1})^{(2n-4)}) \to H^{2n-7}(Z_{n-1}) \to \\
H^{2n-6}_c(\sV(Q_{n-1})^{(2n-4)}\setminus Z_{n-1})\to
  H^{2n-6}(\sV(Q_{n-1})^{(2n-4)}) \to \\
H^{2n-6}(Z_{n-1})
}
Since $H^i(\P^{2n-4}) \surj H^i(\sV(Q_{n-1})^{(2n-4)})$ for $i\le
2n-6$, the lemma follows.
\end{proof}

Now we may put together lemma \ref{lem11.7} and \eqref{11.29} to deduce
\begin{lem}\label{lem11.8} We have
\ga{}{H^{2n-7}(Z_{n-1})(-2) \cong
H^{2n-2}(X_n)/H^{2n-2}(\P^{2n-1});\quad
  n\ge 4 \\
H^2(X_3)/H^2(\P^{5}) \cong H^0(\sV(Q_2)^{(2)})(-2) = \Q(-2). \notag
}
\end{lem}
In order to prove Theorem \ref{thm11.2} it will therefore suffice to prove
\begin{thm}\label{thm11.9} Let
\ml{11.32}{Z_n:= \sV\Big(Q_{n}(B_1,\dotsc,B_{n}, A_1,\dotsc,A_{n-1}), \\
  Q_{n-1}(B_1,\dotsc,B_{n-1},A_1,\dotsc,A_{n-2}),
Q_{n-1}(B_2,\dotsc,B_{n},A_2,\dotsc, A_{n-1})\Big).
}
Then, for $n\ge 3$ we have $H^{2n-5}(Z_n,\Q)) \cong \Q(0).$
\end{thm}
\begin{proof}[Proof of Theorem \ref{thm11.9}] To simplify notation, write
\eq{11.33}{Q_p(i) := Q_p(B_i,\dotsc,B_{i+p-1},A_i,\dotsc,A_{i+p-2}).
}
Given a closed subvariety
$V \subset \P^N$, write $\ell(V) \ge r$ if the restriction maps
$H^i(\P^N) \to H^i(V)$ are surjective for all $i\le r$. (It is equivalent
to require these maps to be an isomorphism for $i\le \min(2\dim V,r)$.) For
$V = \sV(I)$ it is convenient to write $\ell(I) := \ell(\sV(I))$.  For
example a linear subspace has $\ell=\infty$. A disjoint union of
$2$ points has $\ell = -1$.

In what follows, the term {\it variety} is used loosely to mean a reduced
(but not necessarily irreducible) algebraic scheme over a field. We begin
with some elementary properties of $\ell$.
\begin{lem}\label{lem11.10} Let $L \subset \P^N$ be a linear subspace of
dimension
$p$. Let $\pi: \P^N\setminus L \to \P^{N-p-1}$ be the projection with center $L$. For
$V \subset \P^{N-p-1}$ a closed subvariety, write (abusively)
$\pi^{-1}(V) \subset \P^N$ for the cone over $V$. Then
$\ell(\pi^{-1}(V)) = \ell(V)+2(p+1)$.
\end{lem}
\begin{proof}$\pi: \P^N\setminus L \to \P^{N-p-1}$ is an
$\A^{p+1}$-bundle. By homotopy invariance, we have a commutative diagram
\eq{11.34}{\begin{CD}H_c^{i+2(p+1)}(\P^N\setminus L) @>>>
H_c^{i+2(p+1)}(\pi^{-1}(V)\setminus L)
\\ @VV \cong V @VV \cong V \\
H^i(\P^{N-p-1})(-p-1) @>\text{surj.}>> H^i(V)(-p-1).
\end{CD}
}
The bottom horizontal map is surjective for $i\le \ell(V)$,
so the top map is surjective in that range as well. Now consider the
diagram
\eq{11.35}{\begin{CD}0 @>>> H^j_c(\P^N\setminus L) @>>> H^j(\P^N) @>a>> H^j(L) @>>>
0
\\ @. @VV\text{surj.}V @VV c V @| \\
0 @>>> H^j_c(\pi^{-1}(V)\setminus L) @>>> H^j(\pi^{-1}V) @>b>> H^j(L) @>>> 0
\end{CD}
}
Note the maps $a,b$ are surjective in all degrees, so we get short-exact
sequences for all $j$. The left-hand vertical map is surjective if and
only if the central map $c$ is surjective. Since the left hand map is
surjective  for $j \le \ell(V) + 2(p+1)$ by \eqref{11.34}, the lemma
follows.
\end{proof}
\begin{lem}\label{lem11.11} Let $V, W \subset \P^N$ be closed
subvarieties. If
$V\cap W \neq \emptyset$,  then
\eq{11.36}{\ell(V\cup W) \ge \min\Big(\ell(V), \ell(W),2\dim(V\cap
W),\ell(V\cap W)+1\Big). }
\end{lem}
\begin{proof} We use Meyer-Vietoris
\ml{11.37}{ H^{i-1}(V)\oplus H^{i-1}(W) \to H^{i-1}(V\cap
W)\to H^i(V\cup W)
\to \\
H^i(V)\oplus H^i(W) \xrightarrow{g} H^i(V\cap W). }
Note in general if we have $A \subset B
\subset \P^N$, then $H^i(B) \surj H^i(A)$ for $i\le \ell(A)$.
Thus, for $i\le \ell(V\cap W)+1$ we get
\eq{11.38}{0 \to H^i(V\cup W) \to H^i(V)\oplus H^i(W)
\xrightarrow{g} H^i(V\cap W).
}
For $i\le \min(\ell(W), 2\dim(V\cap W))$ the map $g$ above is injective
on $0\oplus H^i(W)$, so $\dim H^i(V\cup W) \le \dim H^i(V)$ and the lemma
follows.
\end{proof}
The proof of Theorem \ref{thm11.9} proceeds by writing
\eq{11.39}{Z_n = \sV(Q_n(1),Q_{n-1}(1)) \cap \sV(Q_n(1),Q_{n-1}(2))
}
from \eqref{11.32}. We remark that the automorphism of projective space
given by
\eq{11.40}{B_1 \mapsto B_n, B_2 \mapsto B_{n-1},\ldots A_1 \mapsto
A_{n-1},
\ldots A_{n-1} \mapsto A_1
}
carries $Q_n(1) \mapsto Q_n(1)$ and $Q_{n-1}(1) \mapsto Q_{n-1}(2)$ so
the varieties on the right in \eqref{11.39} are isomorphic.

\begin{lem}\label{lem11.12} We have
\eq{11.41}{\ell(Q_2(1), Q_{1}(2)) = \ell(Q_2(1), Q_{1}(1)) = \ell(Q_2(1))
= \infty.
}
For $n \ge 3$,
\eq{11.42}{\ell(Q_n(1), Q_{n-1}(2)),\ \ell(Q_n(1), Q_{n-1}(1)),\
\ell(Q_n(1))
\ge 2n-3.
}
\end{lem}
\begin{proof} We write
\eq{11.43}{a_n := \ell(Q_n(1)),\ b_n := \ell(Q_n(1),Q_{n-1}(2))
}
(Using the automorphism \eqref{11.40}, we need
only consider these.). We have
\eq{11.44}{Q_2(1) = B_1B_2-A_1^2,\ Q_1(i) = B_i
}
from which the lemma is immediate in the case $n=2$. For $n=3$ we have the
exact sequence
\eq{11.45}{H_c^i(\P^{3}\setminus \sV(Q_2(2))^{(3)}) \to H^i(\sV(Q_3(1))) \to
H^i(\sV(Q_3(1),Q_2(2)))
}
(cf. \eqref{11.48} below). Since $\ell(\sV(Q_2(2))) = \infty$, the group
on the left vanishes for $i<6$. On the other hand
\ml{11.46}{\sV(Q_3(1),Q_2(2)) = \{B_3=A_2=0\} \cup \{A_1=B_2B_3-A_2^2=0\}
\\
\subset  \P^{4}(B_1,B_2,B_3,A_1,A_2)
}
Each of the two pieces on the right has $\ell = \infty$. Their
intersection is the linear space $L:= \{A_2=A_1=B_3=0\}$ which is a line.
Lemma \ref{lem11.11} gives $b_3:= \ell(Q_3(1),Q_2(2)) \ge 2$, but we can
consider directly the situation for $H^3$
\eq{11.47}{\ldots \surj H^2(L)  \to H^3(\sV(Q_3(1),Q_2(2))) \to 0 \oplus 0
}
and conclude $a_3\ge b_3 \ge 3= \max(3,2\cdot 4-5)$.

The proof of the lemma for $n\ge 4$ is recursive. We have, projecting
from the point $B_1=1, B_i = A_j=0$ using \eqref{11.9},
\begin{small} \minCDarrowwidth.5cm
\eq{11.48}{\begin{CD}H^i_c\Big(\sV(Q_n(1))\setminus \sV(Q_n(1),Q_{n-1}(2))\Big)
@>>> H^i\Big(\sV(Q_n(1))\Big) @>>> H^i\Big(\sV(Q_n(1),Q_{n-1}(2))\Big) \\
@VV\cong V \\
H^i_c(\P^{2n-3} \setminus  \sV(Q_{n-1}(2))^{(2n-3)})
\end{CD}
}
\end{small}
Dropping the variable $A_1$, $\P^{2n-3} \setminus  \sV(Q_{n-1}(2))^{(2n-3)}$ becomes
an $\A^1$-bundle over $\P^{2n-4} \setminus \sV(Q_{n-1}(2))^{(2n-4)}$, so
\eq{11.49}{H^i_c\Big(\sV(Q_n(1))\setminus \sV(Q_n(1),Q_{n-1}(2))\Big) = 0
}
for $i\le a_{n-1}+3$. We conclude from \eqref{11.36} that
\eq{11.50}{a_n \ge \min(a_{n-1}+3,b_n).
}
As a consequence of \eqref{11.9},
\ml{11.51}{(Q_n(1),Q_{n-1}(2)) =
\Big(B_1Q_{n-1}(2)-A_1^2Q_{n-2}(3),Q_{n-1}(2)\Big) = \\
 (A_1^2Q_{n-2}(3),Q_{n-1}(2)).
}
In terms of $\sV$ this reads
\eq{11.52}{\sV(Q_n(1),Q_{n-1}(2)) =\sV(Q_{n-2}(3),Q_{n-1}(2))^{(2n-2)}\cup
\sV(Q_{n-1}(2),A_1)^{(2n-2)}.
}
The varieties on the right are cones with fibres of dimensions $2$ and $1$
respectively. From Lemmas \ref{lem11.10} and \ref{lem11.11} we conclude
\ml{11.53}{b_n \ge \\
\min(b_{n-1} + 4, a_{n-1} + 2,
2\dim\sV(Q_{n-1}(2),Q_{n-2}(3)) + 2, b_{n-1}+3) = \\
\min(a_{n-1}+2,b_{n-1}+3, 4n-10).
}
Starting with $a_3, b_3 \ge 3$ and plugging recursively into \eqref{11.53}
and
\eqref{11.50}, the inequalities of the lemma, $a_n, b_n \ge 2n-3$,
follow.
\end{proof}
We return now to the proof of Theorem \ref{thm11.9}.
\begin{lem}We have the decompositions
\ml{11.54}{\sV(Q_n(1),Q_{n-1}(2)) = \sV(A_1, Q_{n-1}(2)) \cup
\sV(A_2,Q_{n-2}(3))\cup \ldots \\
\cup \sV(A_{n-1},B_n)
}
\ml{11.55}{\sV(Q_n(1),Q_{n-1}(1)) = \sV(A_{n-1}, Q_{n-1}(1)) \cup
\sV(A_{n-2},Q_{n-2}(1))\cup \ldots \\
\cup \sV(A_{1},B_1)
}
\ml{11.56}{\sV(Q_n(1),Q_{n-1}(1)) \cup \sV(Q_n(1),Q_{n-1}(2)) = \\
\sV(A_1,Q_n(1)) \cup \sV(A_2,Q_n(1))\cup\ldots\cup \sV(A_{n-1},Q_n(1)) =
\sV(\prod_{i=1}^{n-1}A_i,Q_n(1)).
}
\end{lem}
\begin{proof} For \eqref{11.54}, we appeal repeatedly to
  \eqref{11.9}
\ml{11.57}{\sV(Q_n(1),Q_{n-1}(2)) = \sV(A_1, Q_{n-1}(2)) \cup
  \sV(Q_{n-1}(2),Q_{n-2}(3)) = \ldots
}
To prove \eqref{11.55}, we apply the automorphism \eqref{11.40} to
\eqref{11.54}. Finally, from the determinant formula \eqref{11.8} one sees
the congruences
\eq{11.58}{Q_n(1) \equiv Q_p(1)\cdot Q_{n-p}(p+1) \mod A_p;\quad 1\le p\le
n-1 }
We can use these to combine the $\sV(A_i,*)$ from \eqref{11.54} and
\eqref{11.55}.
\end{proof}
The idea now is to use Meyer-Vietoris on \eqref{11.39} and \eqref{11.56}.
We get
\ml{11.59}{H^{2n-5}\Big(\sV(Q_n(1),Q_{n-1}(2))\Big) \oplus
H^{2n-5}\Big(\sV(Q_n(1),Q_{n-1}(1))\Big) \to \\
 H^{2n-5}(Z_n) \to H^{2n-4}\Big(\sV(\prod_{i=1}^{n-1}A_i,Q_n(1))\Big)
\to \\
H^{2n-4}\Big(\sV(Q_n(1),Q_{n-1}(2))\Big) \oplus
H^{2n-4}\Big(\sV(Q_n(1),Q_{n-1}(1))\Big) \to H^{2n-4}(Z_n)
}
The vanishing results from Lemma \ref{lem11.12} now yield
\eq{11.60}{H^{2n-5}(Z_n) \cong
H^{2n-4}\Big(\sV(\prod_{i=1}^{n-1}A_i,Q_n(1))\Big)\Big/
H^{2n-4}(\P^{2n-2}).
}
The final step in the proof of Theorem \ref{thm11.9} will be to analyse
the spectral sequence
\eq{11.61}{E_1^{p,q} =
  \bigoplus_{i_0,\dotsc,i_{p}}H^q\Big(\sV(A_{i_0},\dotsc,A_{i_{p}},
  Q_n(1))\Big) \Rightarrow
  H^{p+q}\Big(\sV(\prod_{i=1}^{n-1}A_i,Q_n(1))\Big).
}
We can calculate $H^q\Big(\sV(A_{i_0},\dotsc,A_{i_{p}},Q_n(1))\Big)$ as
follows. Write $n_0 = i_0, n_1=i_1-i_0,\dotsc,n_p = i_p-i_{p-1}, n_{p+1} =
n-i_p$. Thus we have a partition $n = \sum_0^{p+1}n_j$. As in
\eqref{11.58} we may factor
\eq{11.62}{Q_n(1)|_{A_{i_0}=\ldots = A_{i_p}=0} =
Q_{n_0}(1)Q_{n_1}(i_0+1)\cdot
  Q_{n_{p+1}}(i_p+1)|_{A_{i_0}=\ldots = A_{i_p}=0}
}
Each $Q_{n_j}(i_{j-1}+1)$ is a homogeneous function on $\P^{2n_j-2}$. Note
if $n_j=1$, $Q_1(i)=B_i$ is a homogeneous function on $\P^0$. (The
homogeneous coordinate ring of $\P^0$ is a polynomial ring in one
variable.)

We have linear spaces
$$L_j \subset
\P^{2n-p-3}(A_1,\dotsc,\widehat A_{i_0},\dotsc,\widehat
A_{i_p},\dotsc,A_{n-1},B_1,\dotsc,B_n)$$
and cone maps $\pi_j :\P^{2n-p-2} \setminus L_j \to \P^{2n_j-2}$. (When $n_j=1$,
$L_j$ is a hyperplane.) Then
$\sV(A_{i_0},\dotsc,A_{i_p}, Q_n(1))$ is the union of the cones
$\pi_j^{-1}(\sV(Q_{n_j}(i_j)))$. (When $n_j=1$, the cone is just $L_j$.)
Write $$U_j = \P^{2n_j-2} \setminus \sV(Q_{n_j}(i_j))$$ ($U_j = {\rm pt}$ when $n_j=1$)
and $$U = \P^{2n-p-3} \setminus
\bigcup_{j=0}^{p+1} \pi_j^{-1}(\sV(Q_{n_j}(i_j))).$$ The map $\prod \pi_j: U
\to \prod U_j$ is a
$\G_m^{p+1}$-bundle. Thus
\ml{11.63}{H^*_c\Big(\P^{2n-p-3}\setminus \sV(A_{i_0}, \ldots, A_{i_p}, Q_n(1))\Big) = \\
H^*(U) \cong
H^*_c(\G_m^{p+1}) \otimes \bigotimes_{j=0}^{p+1} H^*_c(U_j).
}
Suppose now that some $n_j>1$. Then, by Lemma \ref{lem11.12}, these
cohomology groups vanish in degrees less than or equal to
\eq{11.64}{p+1 + \sum_{j=0}^{p+1}(2n_j-2) = 2n -p -3.
}
It follows that we have surjections
\eq{11.65}{H^i(\P^{2n-2}) \surj
H^i(\sV(A_{i_0},\dotsc,A_{i_p},Q_n(1)));\quad
  i\le 2n-p-4.
}
Note this includes the middle dimensional cohomology.

The exceptional case is when all the $n_j=1$. Then
$p=n-2$. Formula
\eqref{11.64} would suggest $H^*_c(U)=(0), *<n$, but in fact $U \cong
\G_m^{n-1}$ has $H^{n-1}_c(U) \neq 0$. We have
\eq{11.66}{E_1^{n-2,q} = H^q\Big(\sV(A_1,\dotsc,A_{n-1},Q_n(1))\Big) =
H^q(\sV(\prod_{i=1}^n B_i))
}
It follows that $E_2^{n-2,n-2} = \Q$, and $E_2^{p,q}=(0)$ for $p+q=2n-4$,
if $p \neq 0,n-2$. One has
\ml{11.67}{E_2^{0,2n-4} = \ker\Big(\bigoplus_{i=1}^{n-1}
H^{2n-4}(\sV(A_i,Q_n(1))) \to \\
\bigoplus_{I = \{i_1,i_2\}}
H^{2n-4}(\sV(A_{i_1},A_{i_2},Q_n(1)))\Big)
}
Again by \eqref{11.65} $E_2^{0,2n-4}= \Q$ is generated by the class of the
hyperplane section. Finally, the differential $d_r$ reads
\eq{11.68}{E_r^{p-r,q+r-1} \to E_r^{p,q} \to E_r^{p+r,q-r+1}.
}
We have $r\ge 2$.
In the case $p+q=2n-4$, the group on the left vanishes by \eqref{11.65},
the group in the middle vanishes for $p \neq 0, n-2$, and the group on the
right vanishes for $p=n-2$ because we have only $n-1$ components. It
follows that $E_{r+1}^{p,q}\cong E_{r}^{p,q}$. We conclude from
\eqref{11.60}
\eq{11.69}{H^{2n-5}(Z_n) \cong \Q(0).
}
This completes the proof of Theorem
  \ref{thm11.9}.
\end{proof}
By Lemma \ref{lem11.8}, Theorem \ref{thm11.2} follows from Theorem \ref{thm11.9}. This completes the proof of Theorem \ref{thm11.2}.
\end{proof}

\section{de Rham class} \label{sec12}

Let $X_n\subset \P^{2n-1}$ be the graph hypersurface associated to the
wheel and spoke graph with $n$ spokes as in section \ref{sect_ws}. By
the results in that section, we know that de Rham cohomology fulfills $H^{2n-1}_{DR}(\P^{2n-1}\setminus X_n)
\cong K$. Our objective here is to show this is generated by
\eq{1}{\eta_n:=\frac{\Omega_{2n-1}}{\Psi_n^2} \in
\Gamma(\P^{2n-1}, \omega(2X_n))
}
(cf.  \eqref{5.10}), i.e. we show that $[\eta_n] \neq 0$ in
$H^{2n-1}_{DR}(\P^{2n-1}\setminus X_n)$.

To a certain point, the argument is general and applies to the form
$\eta_\Gamma$ attached to any graph with $n$ loops and $2n$ edges. In
this generality it is true that $[\eta_\Gamma]$ lies in the second
level of the coniveau filtration. We do not give the proof here.

\begin{lem} Let $U = \Spec R$ be a smooth, affine variety, and let $0 \neq
f, g\in R$ be functions. Let $Z:f=g=0$ in $U$. We have a map of complexes
\eq{12.2}{\Big(\Omega^*_{R[1/f]}/\Omega^*_R\Big) \oplus
\Big(\Omega^*_{R[1/g]}/\Omega^*_R\Big) \to
\Big(\Omega^*_{R[1/fg]}/\Omega^*_R\Big)
}
Then the de Rham cohomology with supports $H^*_{Z,DR}(U)$ is computed by the
cone of \eqref{12.2} shifted by $-2$.
\end{lem}
\begin{proof}
The localization sequence identifies
\eq{}{H^*_{\{f=0\},DR}(U) = H^*(\Omega^*_{R[1/f]}/\Omega^*_R[-1])
}
(resp. replace $f$ by $g$ resp. $fg$.) The assertion of the lemma follows
from the exact sequence for
$X, Y
\subset U$
\eq{}{\ldots\to H^*_{X\cap Y} \to H^*_X \oplus H^*_Y \to H^*_{X\cup Y} \to
H^{*+1}_{X\cap Y} \to \ldots
}
\end{proof}

\begin{rmk}Evidently, this cone is quasi-isomorphic to the cone of
\eq{}{\Omega^*_{R[1/f]}/\Omega^*_R \to
\Omega^*_{R[1/fg]}/\Omega^*_{R[1/g]}.
}
\end{rmk}
For the application, $U = \P^{2n-1}\setminus X_n$. To facilitate computations, it
is convenient to localize further and invert a homogeneous coordinate as
well. We take $a_i = \frac{A_i}{A_{n-1}}$ and
$b_i=\frac{B_i}{A_{n-1}}$, \eqref{11.4} . (We will check that the forms we work
with have no poles along $A_{n-1} = 0$.)

 We write $Q_p(i)$
as in \eqref{11.33}. Let $q_p(i) = \frac{Q_p(i)}{A_{n-1}^p}$ (resp. $\kappa_n =
\frac{K_n}{A_{n-1}^n}$ with $K_n$ as in \eqref{11.7}). Take $f=q_{n-1}(1),\ g =
q_{n-2}(2)$. The local defining equation $X_n : b_0q_{n-1}(1) + \kappa_n$
has been inverted in $U$, so $\kappa_n$ is invertible on $f=0$ and the
element
\ml{}{\beta:= \\ -db_1\wedge\ldots\wedge db_{n-1}\wedge
da_0\wedge\ldots\wedge
da_{n-2}\frac{1}{\kappa_n}(\frac{1}{q_{n-1}(1)}- \frac{b_0}{b_0q_{n-1}(1)
+ \kappa_n})
}
is defined in $\Omega^{2n-2}_{R[1/f]}/\Omega^{2n-2}_{R}$ and satisfies
\eq{}{d\beta = \eta_n = \frac{db_0\wedge \ldots \wedge db_{n-1}\wedge
da_0\wedge \ldots \wedge da_{n-2}}{(b_0q_{n-1}(1)+\kappa_n)^2}.
}
Applying the fundamental relation expressed by Lemma \ref{lem11.6},
one obtains
\ga{12.7}{\kappa_n q_{n-2}(2)\equiv (a_0q_{n-2}(2) +
(-1)^na_1\cdots a_{n-2})^2 \mod  q_{n-1}(1).}
Computing now in $\Omega^*_{R[1/fg]}/\Omega^*_{R[1/g]}$ we find
\ml{}{\beta = \\
-\frac{dq_{n-1}(1)}{q_{n-1}(1)} \wedge
\frac{db_2}{\kappa_n q_{n-2}(2)}\wedge
db_3\wedge \ldots \wedge da_{n-2} (1-\frac{b_0q_{n-1}(1)}{b_0q_{n-1}(1)
+\kappa_n}) = \\
d\Big(\frac{1}{a_0q_{n-2}(2) +(-1)^n a_1\cdots a_{n-2}} \cdot
\frac{dq_{n-1}(1)}{q_{n-1}(1)}
\wedge\frac{dq_{n-2}(2)}{q_{n-2}(2)} \wedge \nu \Big)
}
where
\ml{10}{\nu=\pm \frac{db_3}{q_{n-3}(3)}\wedge db_4\wedge\ldots \wedge
db_{n-1}\wedge da_1\wedge\ldots\wedge da_{n-2} = \\
\pm \frac{dq_{n-3}(3)}{q_{n-3}(3)} \wedge
\frac{dq_{n-4}(4)}{q_{n-4}(4)}\wedge \ldots
 \wedge \frac{dq_{1}(n-1)}{q_{1}(n-1)}\wedge
da_1\cdots\wedge  da_{n-2}.
}
(Note that $a_0$ is omitted.)

It follows from \eqref{12.7} that in
$\Omega^*_{R[1/fg]}/\Omega^*_{R[1/g]}$ we have
\ml{11}{\beta = d\Big(\frac{1}{a_0q_{n-2}(2) +(-1)^n a_1\cdots a_{n-2}}
\cdot
\frac{dq_{n-1}(1)}{q_{n-1}(1)}
\wedge\frac{db_2}{q_{n-2}(2)} \wedge db_3\ldots \\
 db_{n-1}\wedge
da_1\wedge\ldots\wedge da_{n-2} \Big) = d\theta \\
\theta := \frac{1}{a_0q_{n-2}(2) +(-1)^n a_1\cdots a_{n-2}}
\cdot
\frac{dq_{n-1}(1)}{q_{n-1}(1)}
\wedge\frac{db_2}{q_{n-2}(2)} \wedge db_3\ldots \\
 db_{n-1}\wedge
da_1\wedge\ldots\wedge da_{n-2}
}
(defining $\theta$.) One checks easily that neither $\beta$ nor
$\theta$ has a pole along $A_{n-1}=0$, so the pair
\eq{12a}{(\beta,
\theta) \in H^{2n-1}_{Z, DR}(U)
} represents a class mapping to $\eta_n
\in H^{2n-1}_{DR}(\P^{2n-1}\setminus X_n)$. Here $$Z: Q_{n-1}(1) = Q_{n-2}(2) = 0.$$
\begin{lem}\label{lem3} The map
\eq{12}{H^{2n-1}_{Z}(\P^{2n-1}\setminus X_n) \to
H^{2n-1}(\P^{2n-1}\setminus X_n) }
is injective.
\end{lem}
\begin{proof}Let $Y: Q_{n-1}(1)=0$. We have
\eq{13}{H^{2n-1}_{Z}(\P^{2n-1}\setminus X_n) \xrightarrow{u}
H^{2n-1}_{Y}(\P^{2n-1}\setminus X_n) \xrightarrow{v}
H^{2n-1}(\P^{2n-1}\setminus X_n)
}
and it will suffice to show $u$ and $v$ injective. We have projections
\eq{14}{\P^{2n-1}\setminus (X_n \cup Y) \xrightarrow{B_0} \P^{2n-2}\setminus Y_0
\xrightarrow{A_0, A_{n-1}} \P^{2n-4}\setminus Y_{1}
}
Here $\P^{2n-1}$ has homogeneous coordinates $A_0,\dotsc,A_{n-1},
B_0,\dotsc,B_{n-1}$, the arrows are labeled by the variables which are
dropped, and $Y, Y_0$ are cones over $Y_1$. The arrow on the left is a
$\G_m$-bundle and on the right an $\A^2$-bundle. It follows that
\ml{}{H^{2n-2}(\P^{2n-1}\setminus (X_n \cup Y)) \cong \\
H^{2n-2}(\P^{2n-4}\setminus Y_{1})\oplus
H^{2n-3}(\P^{2n-4}\setminus Y_{1})(-1) = (0)
}
by Artin vanishing. As a consequence, the map $v$ in \eqref{13} is
injective.

The locus $Y\setminus Z$ is smooth ($Q_{n-2}(2) = \partial Q_{n-1}(1)/\partial
B_1$) so to prove injectivity for $u$ it will suffice to show
\eq{}{H^{2n-4}(Y\setminus ((X_n\cap Y)\cup Z)) = (0).
}
Consider the projection obtained as in \eqref{14} by dropping the
variables $B_0, A_0, A_{n-1}$ (so $Y, Z$ are cones over $Y_1, Z_1$)
\eq{}{Y\setminus ((X_n\cap Y)\cup Z) \xrightarrow{\pi} Y_1 \setminus Z_1 \subset
\P^{2n-4} }
Note that $X_n\cap Y:Q_{n-1}(1) = K_n=0$ where $K_n$ is as in
\eqref{11.7}. We can write
$\pi$ as a composition of two projections. First dropping
$B_0$ yields an $\A^1$-fibration. Then dropping $A_0, A_{n-1}$ leads to a
fibration with fibre $\A^2-\text{quadric}$. By Lemma \ref{lem11.6}, this
quadric is a double line, so the fibres of $\pi$ are $\A^2 \times \G_m$.
It follows that
\ml{}{H^{2n-4}(Y\setminus ((X_n\cap Y)\cup Z)) \cong H^{2n-4}(Y_1\setminus Z_1)
\oplus H^{2n-5}(Y_1 \setminus Z_1)(-1) = \\
H^{2n-5}(Y_1 \setminus Z_1)(-1)
}
(The right hand identity is Artin vanishing since $Y_1\setminus Z_1$ is affine of
dimension $2n-5$.) Dropping the variable $B_1$ realizes $\{Q_{n-2}(2)=0\}$
as the cone over a hypersurface $Y_2 \subset \P^{2n-5}$.
Using \eqref{11.9}, we conclude
\eq{19}{H^{2n-5}(Y_1 \setminus Z_1) \cong H^{2n-5}(\P^{2n-5} \setminus Y_2)
}
But the equation defining $Y_2$ does not involve $A_1$, so yet another
projection is possible, and we deduce vanishing on the right in
\eqref{19} by Lemma \ref{lem11.4}.  \end{proof}

\begin{thm} \label{thm12.14} Let $X_n$ be the graph hypersurface for the wheel and spokes
graph with $n$ spokes. Let $[\eta_n] \in H^{2n-1}_{DR}(\P^{2n-1}\setminus X_n)$ be
the de Rham class \eqref{1}. Then $$K[\eta_n]= H^{2n-1}_{DR}(\P^{2n-1}\setminus X_n).$$
\end{thm}
\begin{proof}We have lifted $[\eta_n]$ to a class $(\beta, \theta) \in
H^{2n-1}_{Z,DR}(\P^{2n-1}\setminus X_n)$, \eqref{12a}. By Lemma \ref{lem3}, it
will suffice to show $(\beta, \theta) \neq 0$. We localize at the generic
point of $Z$. It follows from \eqref{10} and \eqref{11} that as a class
in the de Rham cohomology of the function field of $Z$, this class is
represented by the form
\eq{}{\pm d\log(q_{n-3}(3))\wedge\ldots\wedge d\log(q_1(n-1))\wedge
d\log(a_1)\wedge\ldots \wedge d\log(a_{n-2})
}
It is easy to see that this is a non-zero multiple of
$$d\log(b_3)\wedge \ldots \wedge d\log(b_{n-1})\wedge d\log(a_1)\ldots d\log(a_{n-2})$$
and so is nonzero as a form. To see that it is nonzero as a cohomology
class, one applies Deligne's mixed Hodge theory which
implies that the vector space of logarithmic forms injects into de Rham cohomology of the open on which those forms are smooth.
\end{proof}

\section{Wheels and beyond} \label{sec13}
\subsection{A few words on the wheel with $3$ spokes}

Let $X_3 \subset \P^5$ be the hypersurface associated to the wheel with $3$
spokes. $X_3: \det(A_1M_1+\ldots + A_6M_6)=0$ where the $M_i$ are symmetric
rank $1$ $3\times 3$ matrices. It is easy to see in this case that the $M_i$ span
the vector space of all symmetric $3\times 3$-matrices. The mapping
$g\mapsto {}^tgg$ identifies $GL_3(\C)/O_3(\C)$ with the space of
invertible symmetric $3\times 3$ complex matrices. It follows that
\eq{}{\P^5-X_3 \cong GL_3(\C)/\C^\times O_3(\C).
}
From this, standard facts about the cohomology of symmetric spaces yield
theorem \ref{thm11.2} for $X_3$. (We thank P. Deligne for this argument.) 

From another point of view, $X_3$ is the space of singular quadrics in
$\P^2$. Such a quadric is a union of two (possibly coincident) lines, so we
get 
\eq{}{X_3 \cong \text{Sym}^2 \P^2
}
This way we see immediately that
$H^4(X)=\Q(-2)\oplus \Q(-2)$, where the 2 generators are the class of the
algebraic cycles $p\times \P^2 + \P^2 \times p$ and the diagonal
$\Delta$. In particular, Remark \ref{rmk10.5} is clear. 

Then $p\times 
\P^2$ is linearly embedded into $\P^5$ while $\Delta$ is embedded by the
the complete linear system $\sO(-2)$. Thus $\Delta-2 \cdot(p\times \P^2+
\P^2\times p)$
spans the interesting class in $H^4(X)_{{\rm prim}}$.  It is likely that
its strict transform in the blow up $\pi: P\to \P^5$ yields a relative
class in $H^6_Y(P, B)$, but we haven't computed this last piece.

\subsection{Beyond wheels}
An immediate observation is that the wheel with $n$ spokes $w_n$,
\begin{equation} w_n=\wheels\end{equation} and the zig-zag graphs $z_n$, \begin{equation} z_n=\zigzag\end{equation} are both obtained
by gluing triangles together in a rather obvious way. Both classes
of graphs evaluate to rational multiples of $\zeta(2l-3)$ at
$l$-loops \cite{Pisa}. The kinship between these two classes of
graphs is not easily seen at the level of their graph polynomials. Suppose we
try to look directly at the Feynman period \eqref{4.3}. Let $\ell =
e_1+e_2+e_3 \in H_1(\Gamma)$ be the loop spanned by a triangle. If we
choose coordinates on $H_1(\Gamma)$ in such a way that the first coordinate
$k$ coincides on $\Q\cdot\ell \subset H_1$ with $e_i^\vee,\ i\le 3$, and the other
coordinates $q$ are pulled back from a system of coordinates on
$H_1/\Q\cdot\ell$, then the $k$ coordinate appears only in the quadrics
$Q_i$ associated to the edges $e_i,\ i=1,2,3$. Replacing $k$ by
$k_1,\dotsc,k_4$, the period \eqref{4.3} can
be written
\eq{}{\int_{q =-\infty}^\infty \frac{dq}{Q_4(q)\cdots Q_{n}(q)}\int_{k=
    -\infty}^\infty \frac{dk}{Q_1(k,q)Q_2(k,q)Q_3(k,q)}
}

We have the Feynman parametrization
\begin{equation} \frac{1}{Q_1(k)Q_2(k)Q_3(k)}=
\int_0^\infty\!\int_0^\infty
\frac{1+y}{\left[x(1+y)Q_1(k)\!+\!yQ_2(k)\!+\!Q_3(k)\right]^3}dxdy,
\end{equation}
and the elementary integral, valid with appropriate positivity hypotheses
on an inhomogeneous quadric $\widetilde Q(k_1,\dotsc,k_4)$,
\eq{}{\int_{k_1,\dotsc,k_4=-\infty}^\infty \frac{d^4k}{\widetilde Q^3} = \frac{1}{Q}
} 
where, upto a scale factor depending on the determinant of the degree $2$
homogeneous part of $\widetilde Q$, $Q$ is a certain quadratic polynomial in the coefficients of
$\widetilde Q$. With these substitutions, the period becomes
\eq{}{\int_{x,y=0}^\infty dxdy\int_{q=-\infty}^\infty \frac{dq}{Q(x,y,q)Q_4(q)\cdots Q_{n}(q)}
}
where $Q(x,y,q)$ is quadratic in $q$ with coefficients which are rational
functions in the Feynman parameters $x,y$. It would be of interest to try
to make this calculation motivic. 

 A triangle is the one-loop
contribution to the six-point Green function in $\phi^4$ theory:
its four-valent vertices between any pair of its three edges allow
for two external edges, so that these three vertices allow for six
external edges altogether.

 The message in the above that sequences of triangles increase
the transcendental degree (= point at which $\zeta$ is evaluated) in steps of two seems to be a universal
observation judging by computational evidence. Indeed, let us look
at the graph which encapsulates the first appearance of a multiple
zeta value, in this case the first irreducible double sum
$\zeta(5,3)$ which appears in the graph $$M=\Mgraph.$$ This graph
is the first in a series of graphs $$M_i=\Migraphs$$ Adding $\ell$
triangles yields $\zeta(5, 2l+3)$.

Most interestingly, these graphs can be decomposed into zig-zag
graphs in a manner consistent with the Hopf algebra structure on the
multiple zeta value Hopf algebra 
MZVs, upon noticing that the replacement of a triangle in
$$\wthree$$ by the six-point function $$g_6=\bip$$ delivers the graph $M$. (Remove the three edges of a triangle
from $w_3$, and attach the remaining graph, which has $3$ univalent
vertices and one trivalent vertex, to $g_6$ by identifying the univalent vertices with $3$ vertices
of $g_6$ no two of which are connected by a single edge.)
Note that indeed $g_6$ has six vertices of valence three. Each
vertex hence will have one external edge attached to it to make it
four-valent, and the resulting six external edges make this graph
into a contribution to a six-point function. It can hence replace
any triangle. 

 Furthermore, the six-point function $g_6$ is
related to the four-loop graph $$w_4=\wfour$$ by the operation
\begin{equation} w_4=g_6/e\end{equation} where $e$ is any edge
connecting two vertices. Indeed, $g_6$ is the bipartite graph on
two times three edges. Shrinking any of those edges to a point
combines two valence-three vertices into one four-valent vertex
with its four edges connecting to each of the other four remaining
vertices.

This suggests constructing a Hopf algebra $H$ on primitive vertex
graph in $\phi^4$ theory which incorporates the purely
graph-theoretic lemma \ref{lem6.4} such that the following highly
symbolic diagram commutes.
\eq{}{\begin{CD}  H @> \Delta_{\rm 2PI} >> H\otimes H \\
@VV \phi V @VV \phi\otimes \phi V \\
 \rm{MZV} @ > \Delta_{\rm MZV}>> \rm{MZV} \otimes \rm{MZV}
\end{CD}}
First results are in agreement with the expectation that all
graphs up to twelve edges are mixed Tate, which they are by
explicit calculation \cite{Pisa}, and also predict correctly the
apperance of a double sum $\zeta(3,5)$ or products
$\zeta(3)\zeta(5)$ in six-loop graphs. The seven loop data demand
some highly non-trivial checks (currently in process) on the data amassed in
\cite{Pisa,BK}.

\newpage
\bibliographystyle{plain}
\renewcommand\refname{References}

\end{document}